\documentclass[12pt]{amsart}

\usepackage[colorlinks=true,urlcolor=blue,
bookmarks=true,bookmarksopen=true,citecolor=blue
]{hyperref}

\usepackage{pdfsync}
\usepackage[all, cmtip]{xy}
\usepackage{graphicx}

\usepackage{epsfig}
\usepackage{amscd}
\usepackage[mathscr]{eucal}
\usepackage{amssymb}
\usepackage{amsxtra}
\usepackage{amsmath}
\usepackage{latexsym} 
\usepackage[all]{xy}
\usepackage{enumerate}
\usepackage{mathrsfs}

\usepackage{color}
\theoremstyle{plain}

\newtheorem{thm}{Theorem}[section]

\newtheorem{prop}[thm]{Proposition}

\theoremstyle{definition}
\newtheorem{dfn}[thm]{Definition}

\theoremstyle{remark}
\newtheorem{rem}[thm]{Remark}

\theoremstyle{plain}

\newcommand{\cobto}{\leadsto}

\newcommand{\R}{\mathbb{R}}
\newcommand{\Z}{\mathbb{Z}}

\newcommand{\C}{\mathbb{C}}

\newcommand{\mor}{{\textnormal{Mor\/}}}

\DeclareMathOperator{\Ham}{Ham}
\DeclareMathOperator{\osc}{osc}

\newcommand{\ocaddress}{cornea@dms.umontreal.ca}

\newcommand{\egaddress}{egorshel@ias.edu}

\begin{document}

\title[Lagrangian cobordism and metric invariants]{Lagrangian cobordism and metric invariants.}
\date{\today}

\thanks{The first author was supported by an NSERC Discovery grant, a FQRNT Group Research grant,
a Simons Fellowship and an Institute for Advanced Study fellowship grant. The second author was partially supported by the National Science Foundation under agreement No. DMS-1128155.  }

\author{Octav Cornea and Egor Shelukhin}

 \address{Octav Cornea, Department of Mathematics
  and Statistics, University of Montreal, 
  Montreal, QC H3C 3J7 and Institute for Advanced Study, Einstein Drive, Princeton NJ 08540} \email{\ocaddress}

\address{Egor Shelukhin, Institute for Advanced Study, Einstein Drive, Princeton NJ 08540} \email{\egaddress}

\bibliographystyle{plain}


%

\begin{abstract}
We introduce new pseudo-metrics on spaces of Lagrangian submanifolds of a symplectic manifold $(M,\omega)$ by considering areas associated to projecting Lagrangian cobordisms in $\C \times M$ to the ``time-energy plane" $\C$. We investigate the non-degeneracy properties of these pseudo-metrics, reflecting the rigidity and flexibility aspects of Lagrangian cobordisms. 
\end{abstract}

\maketitle

%
%

\tableofcontents 


\section{Introduction}



One of the central aims of symplectic topology is to understand the topology
of the Lagrangian submanifolds  $L$ of a given symplectic manifold $(M^{2n}, \omega)$ which is closed or tame at infinity. A  basic question is whether there is some natural topology, or a distance, on 
the space $\mathcal{L}(M)$ of all such submanifolds that has some interesting, specifically 
symplectic, features.  Without additional constraints on the class of Lagrangians under consideration,
a positive answer to this question is hard to expect for two main reasons: the topological type of the submanifolds $L$ in question is not fixed; symplectic rigidity properties are not preserved by isotopies, even Lagrangian ones,
but only by Hamiltonian isotopies. 

\

Recall that the Hamiltonian diffeomorphisms group, $\Ham(M,\omega)$,  is endowed with a bi-invariant metric 
introduced by Hofer. Many recent advances in symplectic topology are related to properties of Hofer's geometry.  The Hofer metric is also relevant to our problem: if we fix some $L\in \mathcal{L}(M)$ and consider the orbit $\mathcal{L}^{L}(M)$ of $L$ under the action of $\Ham(M,\omega)$ (in other words, these are all the Lagrangians in $M$ that are Hamiltonian isotopic to $L$), then the Hofer metric adjusts naturally to this context and, as noted by Chekanov \cite{ChekanovFinsler,Chekanov:Lag-energy-1} , provides a metric on  $\mathcal{L}^{L}(M)$.

\

The purpose of this paper is to show that there is a natural construction of a family of metrics 
that are defined and non-degenerate on certain subsets of $\mathcal{L}(M)$ that are, in general, larger than Hamiltonian orbits. Moreover, this construction extends the construction of the Hofer 
metric. We will see that using such a metric one can compare certain Lagrangians that are not even
smoothly isotopic.

\subsection{Measuring cobordisms} \label{subsec:cob}
The construction of the metrics mentioned above employs the Lagrangian cobordism
machinery as developed in \cite{Bi-Co:cob1,Bi-Co:cob2}. The notion of Lagrangian cobordism 
was introduced by Arnold \cite{Ar:cob-1,Ar:cob-2} and we refer to \cite{Bi-Co:cob1} for the variant that we use in this paper. In short, a Lagrangian cobordism $V$ is a non-compact Lagrangian submanifold of $(\C\times M,\Omega=\omega_{0}\oplus \omega),\; \omega_0 = dx \wedge dy$ whose ends are {\em cylindrical} of two types, positive and negative.  The
 positive ones are of the form $[\beta_{+}, \infty)\times \{r\}\times L'_{r}$ where $r\in \Z_{>0}$, $\beta_{+}\in \R$  and $L'_{r}$ is a Lagrangian in $M$. Similarly, the negative ends are of the form $(-\infty, \beta_{-}]\times \{r\}\times L_{r}$ with $\beta_{-}<\beta_{+}$. 
 We refer to each of the $L'_{j}$'s as a positive end of $V$ and to each of the $L_{i}$'s as a negative end. A cobordism like this
is sometimes written as  $V: (L'_{j})_{1\leq j\leq k_{+}}\cobto (L_{i})_{1\leq i\leq k_{-}}$ when
we have $k_{+}$ positive ends and $k_{-}$ negative ends. In case $k_{+}=k_{-}=1$ the cobordism is
called {\em simple}. 

\begin{figure}[htbp]
   \begin{center}
      \epsfig{file=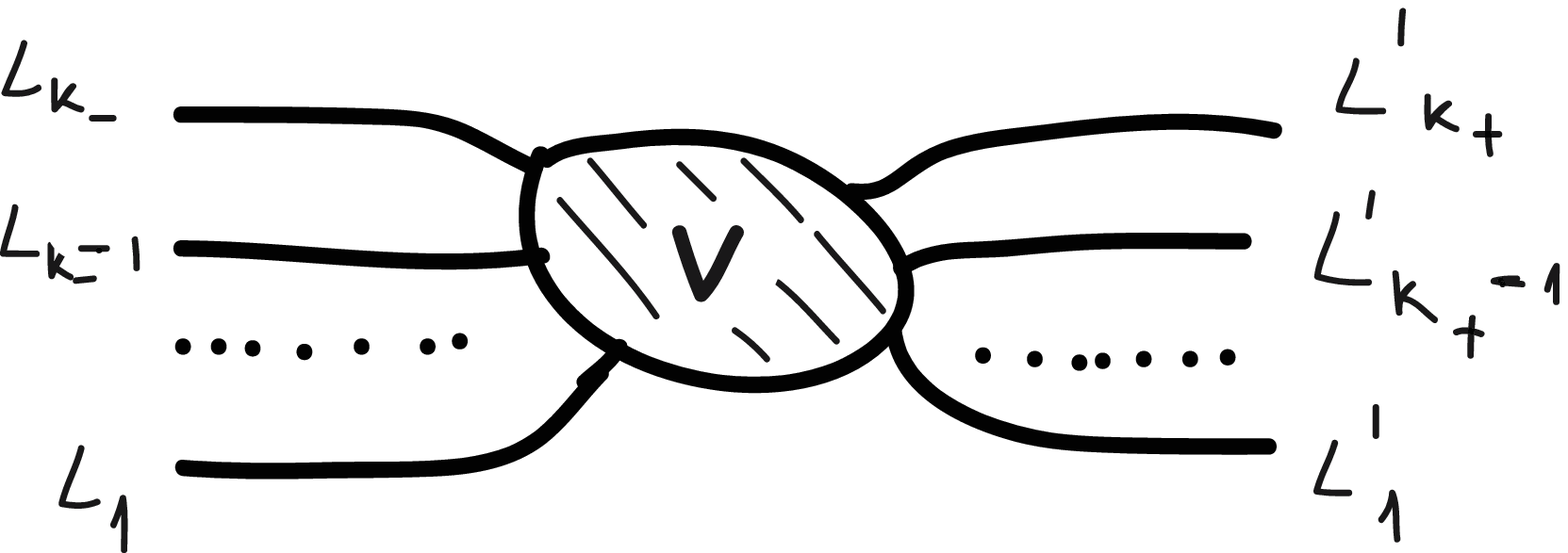, width=0.6\linewidth}
   \end{center}
   \caption{\label{fig:Cobd} A cobordism $V:(L'_{j})\cobto (L_{i})$ projected
     on $\mathbb{C}$.}
\end{figure}

Under the canonical projection $\pi:\C\times M\to \C$ a cobordism looks as in Figure \ref{fig:Cobd}.

Our construction is based on the following
natural  measure of Lagrangian cobordisms.

\begin{dfn}\label{def:shadow} Given a Lagrangian cobordism $V\subset \C\times M$, $V:(L'_{j}) \cobto (L_{i})$, the {\em outline} of $V$, $ou_{V}\subset \C$,  is the closed subset of $\C$ given as the complement of the union of the unbounded components of $\C\backslash \pi(V).$ The {\em shadow} of $V$ is given by:
$$\mathcal{S}(V)=\mathrm{Area\/}(ou_{V})~.~$$
\end{dfn}


\subsection{Statement of the main result}
There is a  dichotomy
\begin{center} {\em flexibility} $\ \leftrightarrow\ $ {\em rigidity}\end{center} in symplectic topology that is very much in evidence in the study of Lagrangian submanifolds and which is also apparent in the main result of the paper, Theorem \ref{thm:metrics} below. This dichotomy is reflected in the way certain topological constraints impact the geometric (symplectic) properties of Lagrangian submanifolds. We work in this paper with certain classes  $\mathcal{L}^{\ast}(M)$ of Lagrangian submanifolds in $(M,\omega)$ as well as with the corresponding classes of simple cobordisms $\mathcal{L}^{\ast}_{cob}(\C\times M)$.  The superscript $(-)^{\ast}$ refers to the constraints on the Lagrangians considered. These constraints are essentially
standard in the field but we refer to \S\ref{subsec:class-Lagr} for details.  In short, from flexible to rigid, $\ast$ can take the following values: $g$ for general, or unconstrained Lagrangians and cobordisms; $w-m$ for weakly monotone by which we mean monotone without restrictions on 
the minimal Maslov class (with a fixed monotonicity factor); $m$ for monotone; $e$ for exact; finally, $L_{0}$ for some fixed Lagrangian $L_{0}\subset M$ with the notation $\mathcal{L}^{L_{0}}(M)$ meaning the Hamiltonian orbit of 
$L_{0}$ (with cobordisms given by Lagrangian suspensions) as in the beginning of the introduction. An inequality such as $\ast\geq w-m$ means that we work with Lagrangians that are at least as rigid as weakly-monotone.

\begin{thm}\label{thm:metrics} 
For two Lagrangians $L,L'\in \mathcal{L}^{\ast}(M)$ let:
 \begin{equation}\label{eq:distances}
d^{\ast}(L,L')=\inf \ \{ \mathcal{S}(V) \ | \ V:L\cobto L' \ , \ V\in \mathcal{L}^{\ast}_{cob}(\C\times M) \}~.~
\end{equation}
\begin{itemize}
\item[i.] If $\ast\geq w-m $, then
$d^{\ast}(L,L')$ defines a metric, possibly infinite, on $\mathcal{L}^{\ast}(M)$. 
\item[ii.] Given any Lagrangian $L_{0}\subset M$ the metric $d^{L_{0}}$ on the Hamiltonian orbit of $L_{0}$ is the Lagrangian analogue of  Hofer's metric.
\item[iii.] There are Lagrangians $L,L'\in \mathcal{L}^{w-m}(M)$ so that $d^{w-m}(L,L')\not=\infty$,
$L$ and $L'$ are not smoothly isotopic and $d^{m}(L,L')=\infty$.
\item[iv.] For $\ast=g$, equation (\ref{eq:distances}) defines a degenerate pseudo-metric. More precisely, 
for any two Lagrangians $L,L'$ the only possible values of $d^{g}(L,L')$ are $0$ and $\infty$ and
$d^{g}(L,L')=\infty$ iff there are no simple Lagrangian cobordisms relating $L$ to $L'$. 
\end{itemize}
\end{thm}
We emphasize that $d^{\ast}(L,L')=\infty$ means that $L$ and $L'$ are not
cobordant via a simple cobordism of class $\ast$ and $d^{\ast}(L,L')=0$ means 
that there are simple cobordisms of arbitrarily small shadow that belong to this class and relate $L$ and $L'$. 

\

The structure of the paper is as follows. Section \ref{subsec:class-Lagr} contains the required background
on the various constraints $\ast$. In \S\ref{sec:tech} we state and prove the main technical ingredient needed
to establish Theorem \ref{thm:metrics}. This result is of some independent interest by itself. Its 
 application in this paper is in \S\ref{sec:proofs} where we show Theorem \ref{thm:metrics}. The paper
 ends with a few other comments in \S\ref{sec:comm}.
 
 \subsection*{Acknowledgements} This paper has grown from a number of questions that have emerged 
 in discussions between Paul Biran and the first author during their long-term investigation of Lagrangian cobordism. We thank him as well as Misha Entov, Kenji Fukaya, and Misha Khanevsky for useful discussions.
 We thank Emmy Murphy for suggesting that the pseudo-metric $d^{g}$  is degenerate by a simple application of the h-principle.

\section{Background: classes of Lagrangians, from flexible to rigid} \label{subsec:class-Lagr}
In this section we recall a series of standard, more and more strict topological constraints on Lagrangian manifolds and fix the relevant notation.
 
\subsection{ $\ast = g$} \label{subsubsec:gen} This is the most flexible choice of contraint and indicates that no restriction is 
imposed - $g$ comes from {\em general}. In other words, $\mathcal{L}^{g}(M)$ are all the Lagrangian submanifolds of $M$ and $\mathcal{L}^{g}_{cob}(\C\times M)$ are all the simple cobordisms in $\C\times M$. 

\subsection{$\ast =w-m$} \label{subsubsec:wmo} This is the {\em weakly - monotone} case. Given a Lagrangian $K\subset M$
there are two morphisms:
$$\mu: \pi_{2}(M,K)\to \Z \ \ \ \ \ \omega : \pi_{2}(M,K)\to \R$$
the first being given by the Maslov class and the second by the integration of $\omega$.
The Lagrangian $K$ is called weakly-monotone if there exists a constant $\rho\in \R$ so that for each
class $\alpha\in \pi_{2}(M,K)$ we have $\omega(\alpha)=\rho \mu (\alpha)$. Notice that there
are no restrictions in this case on the constant $\rho$ or on the minimal Maslov class. The Lagrangians
in the class $\mathcal{L}^{w-m}(M)$ are all weakly monotone with the same monotonicity constant $\rho$
and similarly for the cobordisms in the class $\mathcal{L}^{w-m}_{cob}(\C\times M)$. To be more precise
we may include the constant $\rho$ in the notation in which case we write $\mathcal{L}^{w-m(\rho)}(M)$.

\subsection{$\ast=m$} \label{subsubsec:mo}This is the {\em monotone} case.  A Lagrangian $K\subset M$ is called monotone if it is weakly-monotone and, additionally, the constant $\rho \geq 0$ and, further, the minimal
Maslov class:
$$N_{K}= \min\{ \mu(\alpha) \ | \ \alpha\in \pi_{2}(M,K) \ , \ \omega(\alpha)>0\}$$
is at least $2$. In this context Floer homology is well-defined \cite{Oh:HF1,Oh:HF1-add}. For simplicity we will
mainly deal in this paper with Floer invariants over the base ring $\Z_{2},$ which allows us to disregard questions of orientation.
 If $K$ is monotone, then for a point $P\in K$ and a generic almost complex structure $J$ on $M$,
the number of $J$-holomorphic disks going through $P$ - these are the maps $u:D^2 \to M$, $u(S^{1})\subset K$, $u(+1)=P$, $\bar{\partial}_{J} u=0$ modulo reparametrization - is finite. We denote by
$d_{K}\in \Z_{2}$ the number of these disks modulo $2$. The set $\mathcal{L}^{m}(M)$ indicates
in this case all the monotone Lagrangians in $M$ with the same monotonicity constant $\rho$ as
well as with the same number $d_{K}$. In case we want to indicate explicitly these two constants
we will write $\mathcal{L}^{m(\rho,d)}(M)$ to mean the monotone Lagrangians $K$ in $M$ of monotonicity constant $\rho$ and so that $d_{K}=d$. 
We use similar notation for cobordisms so that, for instance, $\mathcal{L}^{m(\rho, d)}_{cob}(\C\times M)$ are the cobordisms in $\C\times M$ that are 
monotone as Lagrangians in $\C\times M$ with the respective constants $(\rho, d)$.
For any two Lagrangians in $\mathcal{L}^{m(\rho,d)}(M)$  the Floer homology of the pair is well-defined \cite{Oh:HF1,Oh:HF1-add}. To fix some relevant notation we first indicate that unless mentioned otherwise we work at all times over $\Z_{2}$.
For any two  $L,L'\in \mathcal{L}^{m(\rho,d)}$ the Floer chain complex of $L$ and $L'$ will be denoted by $CF(L,L')$. The homology of this complex, $HF(L,L')$, is called the Floer 
homology of $L$ and $L'$. If $L$ and $L'$ intersect transversely, then
$CF(L,L')$ is basically identified with the $\Z_{2}$ vector space generated by the intersection points of $L$ and $L'$ and the differential of the Floer complex counts (with weights given by the symplectic area) $J$-holomorphic strips joining these intersection points.  A more precise description appears in \S\ref{sec:tech}. Floer homology can also be adapted to the case of cobordisms themselves - as in \cite{Bi-Co:cob1,Bi-Co:cob2} -
so that, given any two cobordisms $V, V'\in \mathcal{L}_{cob}^{m(\rho,d)}(\C\times M)$
the Floer complex $CF(V,V')$ is well-defined (as usual, up to canonical quasi-isomorphisms).
There is also an obvious notion of {\em Lagrangian with cylindrical ends} in $\C\times M$  that is more 
general than cobordism. This is defined in the same way as in the cobordism case
except that the ends are of the type $(-\infty, \beta]\times \{a_{i}\}\times L_{i}$, respectively $[\alpha,\infty)\times \{b_{i}\}\times L'_{j}$ for $a_{i},b_{j}\in\R$ while for cobordisms $a_{i},b_{j}\in \Z_{>0}$.
Floer homology is again defined for any pair of such Lagrangians. Furthermore,
there is also an associated notion of isotopy for cobordisms \cite{Bi-Co:cob1} (as well as, more generally,
for Lagrangians with cylindrical ends): two cobordisms $V,V' \subset \C \times M$ are {\em horizontally isotopic} if there exists a Hamiltonian isotopy $\{\phi_{t}\}_{t \in [0, 1]}$ of $\C \times M$ sending $V$ to $V'$ and so that, outside of a compact, $\phi_{t}(V)$ has the same ends as $V$ for all $t \in [0, 1]$ (in other words, the ends can slide along but their image in $\C \times M$ - outside a large compact set - remains the same; the Hamiltonian isotopy is not necessarily with compact support).  As shown in \cite{Bi-Co:cob1}, this type of Hamiltonian
isotopy leaves invariant the Floer homology $HF(V,V')$ just as in the usual compactly supported setting.
The distinction between cobordisms and Lagrangians with cylindrical ends seems somewhat arbitrary but 
is relevant, in fact, for the definition of the Fukaya category of Lagrangian cobordisms - such as in \cite{Bi-Co:cob2}. As this Fukaya category is not needed in this paper, the two notions will be used interchangeably.

\subsection{$\ast=e$}\label{subsubsec:ex} This is the {\em exact} case.

In this case, we assume that the manifold $M$ is exact, i.e. $\omega = d\lambda,$ and moreover, the primitive $\lambda$ restricts to an exact form on each of the Lagrangians belonging to $\mathcal{L}^{e}(M),$ and similarly for cobordisms. The Floer machinery 
was initially developed in the exact case \cite{Fl:Morse-theory} and, as mentioned above, this was later extended
to the monotone setting.
 \subsection{$\ast=L_{0}$}\label{subsubsec:orb} This is the case of Hamiltonian orbits. We fix $L_{0}$ a 
Lagrangian in $M$ and we denote by $\mathcal{L}^{L_{0}}(M)$ all the Lagrangians in $M$ that are Hamiltonian isotopic to $L_{0}$.  The cobordisms in $\mathcal{L}_{cob}^{L_{0}}(\C\times M)$
consist only of Lagrangian suspensions.

We now recall the Lagrangian suspension construction. 
Let $\phi \in \Ham(M,\omega)$ be a Hamiltonian diffeomorphism.  Let $G:[0,1]\times M\to \R$ be a time dependent Hamiltonian so that the time-$1$ diffeomorphism associated to $G$ is $\phi$, $\phi^{G}_{1}=\phi$. We denote by $\phi^{G}_{t}$ the time-$t$ Hamiltonian diffeomorphism associated to $G$ for all $t \in [0,1]$. In particular $\phi^{G}_{0}=id$.  We may assume that $G$ is normalized in such a way so that for some small $\epsilon >0$ it vanishes on $([0,\epsilon]\cup [1-\epsilon, 1])\times M$. 
Such a normalization is easy to achieve by reparametrizing the Hamiltonian flow (see \cite{Po:hambook}):  if $b:[0,1]\to [0,1]$ is a smooth function, the  Hamiltonian isotopy $\phi^{G}_{b(t)}$ is
generated by the Hamiltonian function $(t,z)\to b'(t)G(b(t), z)$ and we may take
$b(t)$ equal to $t$ on $[\frac{3\epsilon}{2}, 1-\frac{3\epsilon}{2}]$ and constant, equal to $0$ on $[0,\epsilon]$ and constant equal to $1$ on $[1-\epsilon, 1]$. In view
of our normalization,  we extend $G$ by zero outside of $[0,1]\times M \subset \R \times M,$ and view it
as a map $G:\R\times M\to \R$.  There is a symplectomorphism $\Phi:  \C\times M\to \C\times M$ defined by $$\Phi(x+iy,z) = (x + i(y+G(x,\phi^{G}_{x}(z))), \phi^{G}_{x}(z))~.~$$
This symplectomorphism is itself Hamiltonian (but not horizontal). For convenience we will denote the 
corresponding Hamiltonian isotopy by $\Psi_{t}:\C\times M\to \C\times M$ so that $\Psi_{1}=\Phi $
and $\Psi_{0}=id$. It is also possible to assume (by an appropriate reparametrization) that the path of
symplectomorphisms $\Psi_{t}$ is constant for $t$ near the ends of the interval $[0,1]$.
Fix now a connected, closed Lagrangian $L\subset M$ and let its Lagrangian suspension along $G$ (see \cite{Po:hambook} as well as  \cite{Cha-Co:SeRep} where the cobordism perspective on this construction is made explicit) be defined by:
$$L^{G}=\Phi(\R\times L)~.~$$ 
Our normalization for $G$ implies that $L^{G}$ is a Lagrangian cobordism, $L^{G}:\phi (L)\cobto L$. Define $\mathcal{L}_{cob}^{L_{0}}(\C\times M)$ to be the set of Lagrangian cobordisms in $\C \times M$ that can be written as a Lagrangian suspension $L^{G}$ for some $L\in \mathcal{L}^{L_{0}}(M)$ and $G:[0,1]\times M\to \R$.

In case $L_{0}$ is monotone, then obviously $\mathcal{L}^{L_{0}}(M)\subset \mathcal{L}^{m(\rho,d)}(M)$ where $\rho$ is the monotonicity constant of $L$ and $d=d_{L}$ and, in particular, all the 
Floer machinery applies. However, if $L_{0}$ is not restricted in any way,
the Floer homology $HF(L_{0}, L_{1})$ is not defined in general even if $L_{1}$ is Hamiltonian isotopic
to $L_{0}$, that is $L_{1}=\phi (L_{0})$ for some $\phi\in \Ham(M,\omega)$.  Despite this, 
starting with Chekanov's work \cite{Chekanov:Lag-energy-1}, it is now well-known that at least parts of the Floer machinery can be used inside $\mathcal{L}^{L_{0}}(M)$ for energies under the bubbling threshold.  In all cases, from the perspective of this paper,
the restriction to just one Hamiltonian orbit is the most rigid constraint in our list.

\

Given that we have listed the possible choices for $\ast$ in order, starting from the most flexible to the most rigid we will say, for instance, that a Lagrangian $L$, or a family of Lagrangians, is at least weakly-monotone to mean that they belong to $\mathcal{L}^{\ast}(M)$ for some $\ast\geq w-m$ - in this case the value of 
 the constraint $\ast$ could be $w-m,m,e$ or $L_{0}$ -  and similarly for cobordisms.

\section{The main technical result}\label{sec:tech}
The proof of the main Theorem is based on a technical  result whose statement and proof are contained in this section. To formulate it we need to recall two additional numerical invariants associated to Lagrangian submanifolds.

\

The first is a positive real number 
associated to a pair of two Lagrangian submanifolds $L, L'\subset M$. It is called
the {\em Gromov width of $L$ relative to $L'$} (while we use this terminology for the sake of brevity, a more appropriate name could be ``the Gromov width of $L$ inside $M\backslash (L'\backslash L)$'') and was introduced in \cite{Bar-Cor:Serre}:
\begin{eqnarray}\label{eq:Gr-w}
w(L,L')=\sup_{r>0}\ \ \{\frac{\pi r^{2}}{2} \ | \ \exists \ e: B_{r}\to M,   \ \mathrm{symplectic\ embedding,}\hspace{0.9in}\\ \nonumber \hspace{0.9in} \ e^{-1}(L)=B_{r}\cap \R^{n}, \ e(B_{r})\cap L'=\emptyset \}~.~
\end{eqnarray}
Here $B_{r}\subset \C^{n}$ is the standard ball of radius $r$ and center $0$.

\

The second is a positive number associated to a Lagrangian $L$ and an almost 
complex structure $J$ on $M$ that is compatible with $\omega$. It is 
called the {\em bubbling threshold of $L$ with respect to $J$}:
\begin{equation} \label{eq:threshold}
\delta(L;J)=\inf_{u}\ \{ \omega(u) \ | \  u \ \mathrm{is\ J-holomorphic}, u:(D^{2}, S^{2})\to (M,L) \  \mathrm{or}\ u: S^{2} \to M\}~.~
\end{equation}
Obviously, this represents the energy at which the bubbling of either a $J$-disk with boundary on $L$ or of a $J$-sphere may occur. It is taken to be $\infty$ if there are no relevant $J$-disks or spheres. 
This definition also make sense in case the almost complex structure $J$ is time-dependent, $J=\{J_{t}\}_{t\in[0,1]}$. In this case, we take the infimum in Equation (\ref{eq:threshold}) over all
disks $u:(D^{2}, S^{1})\to (M,L)$ that are $J_{1}$-holomorphic and over all spheres $u:S^{2}\to M$
that are $J_{t}$-holomorphic for some $t\in [0,1]$. This is a bit of an ad hoc notion - as we do not take into
account the $J_{0}$-holomorphic disks with boundary on $L$ - but it is convenient
for our purposes here (and it provides a better lower bound). 

This number is well-defined also for cobordisms $V$ in case the almost complex
structure $J$ on $\C\times M$ is compatible with $\Omega$ and has the following properties: 
\begin{itemize}
\item[i.]there is a compact family of almost complex structures $\mathcal{J}$ on $M$ compatible with
$\omega$ and a set $K\subset \C$ 
so that for  $z\in \C\setminus K$, $J$ is of the form $i\times J'$ with $J'\in \mathcal{J}$. 
\item[ii.] for some $\alpha_{-}>0$ it restricts to $i\times J_{-}$, $J_{-}\in \mathcal{J}$,  over the set 
$(-\infty,-\alpha_{-}]\times \R\times M$.
\item[iii.] for some $\alpha_{+}>0$, $J$ restricts to $i\times J_{+}$, $J_{+}\in\mathcal{J}$
 over the set $[\alpha_{+}, \infty)\times \R\times M$.  
 \end{itemize}
For convenience, we will
also assume (without loss of generality) that the compact set $K$ above is included in $[-\alpha_{-},\alpha_{+}]\times\R$. We call such almost complex structures trivial at infinity.
We will say that $J_{-}$ is the negative end of $J$ and that $J_{+}$ is the positive end of $J$.

\begin{prop}\label{prop:shadow-est}
Let $V:(L_{1},\ldots, L_{i-1},L, L_{i}\ldots, L_{r})\cobto (L'_{1},\ldots, L'_{j-1},L', L'_{j}\ldots, L'_{s})$ be a cobordism.  
If $L\cap L_{k}=\emptyset$\  \ $\forall\ 1\leq k\leq r$ and $L\cap L'_{m}=\emptyset$\  \ $\forall\ 1\leq m\leq s$, then for each $\epsilon_0 > 0$ there exists a time independent almost complex
structure $J_{-}$ on $M$ (depending on $\epsilon_0$) that is compatible with $\omega$ so that
$$\mathcal{S}(V)\geq \min\ \{ w(L,L') - \epsilon_0 , \delta(V;J)\}$$
 for any time dependent almost complex structure $J= \{J_{t}\}_{t\in [0,1]}$ on $\C\times M$
with the following three properties:
\begin{itemize}
\item[a.] $J$ is compatible with $\omega_{0}\oplus \omega$,
\item[b.] $J_{t}$ is trivial at infinity with $J_{-}$ as its negative end for all $t\in [0,1]$,
\item[c.] $J_{0}=i\times J_{-}$. 
\end{itemize}
\end{prop}

\begin{rem} \label{rem:tech}
a. One reason why this result is of interest is that there are no conditions of any sort imposed
on the Lagrangians and cobordisms involved (or, in the terminology of the paper, $\ast=g$). Thus, 
in practice, to ``measure'' (estimate from below) the distance between $L$ and $L'$ using the shadow of cobordisms the whole question comes down to having a uniform lower bound for $\delta(V,J)$ that applies to all the cobordisms
$V$ in a given class.

b. The presence of $\epsilon_{0}$ in the inequality claimed in the Proposition is required  for the following reason. The almost complex structure $J_{-}$  basically depends on the choice of an  embedding $e:B_{r}\hookrightarrow M$, as in the definition of relative width, so that
 $\frac{\pi r^{2}}{2}\geq w(L,L')-\frac{\epsilon_{0}}{2}$. With this choice, $J_{-}$ extends the standard almost complex structure $e_{\ast}J_{0}$ outside $B_{r}$. In particular, the number $\delta(V,J)$ also depends on the  choice of $e$.  Certainly, by picking different embeddings $e$ for increasing values of $r$ we may reduce
 $\epsilon_{0}$ arbitrarily close to $0$.  However, in this process the choices of $J_{-}$ vary and to 
 eliminate $\epsilon_{0}$ from the statement we would need to control the convergence
of the associated $\delta(V,J)$'s. While this might be possible, we prefered to avoid this additional complication
 as it is not justified in view of our applications here.

c. The statement of the proposition specialized to the case of simple cobordisms (i.e. $r=0=s$), is the only case needed to prove Theorem \ref{thm:metrics}. The non-simple case is in fact an immediate
consequence of the proof in the simple case.
Additionally, the conditions on the ends of the cobordism $V$ (different from $L$ and $L'$) are quite stringent and it is not clear whether these precise conditions are necessary. However, it is certain is that some conditions that ``separate'' $L$ and $L'$ from the other ends
need to be imposed.  To see this consider two Lagrangians $L,L'\subset M$ that are disjoint. 
We can easily find such Lagrangians, even exact, in certain symplectic manifolds. One example is provided by two
homologically non-trivial, disjoint curves on a surface of genus $2$ (note that these Lagrangians are weakly exact).
Consider two curves $\gamma, \gamma'\subset \C$ so that: $\gamma=\R+2i$; $\gamma'$ intersects $\gamma$ in the single point $2i\in\C$; outside of $[-10,10]\times \R$ $\gamma'$ 
equals  $\{[10,+\infty)+3i \}\cup \{(-\infty, -10]+i\}~.~$  
We now consider the Lagrangian $W=(\gamma\times L) \cup (\gamma'\times L')$. As $L,L'$ are disjoint this 
is a cobordism $W:(L',L)\cobto (L,L')$ with vanishing shadow. This example $W$ is not connected but by using $L$ and $L'$ that intersect transversely at a single point (such as the longitude and latitude on a torus) we can start with the Lagrangian $W$ constructed as before - which is now immersed - use Lagrangian surgery  (\cite{La-Si:Lag},\cite{Po:surgery}) to eliminate the single self-intersection point of $W$ and, by using a sufficiently small Lagrangian handle in the surgery,
obtain for any $\epsilon>0$ a connected cobordism $V_{\epsilon}:(L',L)\cobto (L,L')$ of shadow $\leq \epsilon$. 

c. For certain choices of $\ast$, simple cobordisms verifying the constraint $\ast$ are quite special. 
For instance, it is conjectured that in the exact case, $\ast=e$, any simple cobordism is horizontally Hamiltonian isotopic to a Lagrangian suspension.  A partial result in this direction is due to
 Suarez-Lopez \cite{Sua:thesis} who showed that under some topological constraints any exact simple cobordism is smoothly trivial. In the monotone case, $\ast=m$,  only very recently there has been constructed - by Haug \cite{Haug} - a simple monotone cobordism (with Maslov class at least $2$) with non-homeomorphic ends. It is useful to note that if two Lagrangians
 $L, L'\in\mathcal{L}^{m}(M)$ are related by a simple cobordism in $\mathcal{L}^{m}_{cob}(\C\times M)$, then, by the results in \cite{Bi-Co:cob2}, $L$ and $L'$ are isomorphic objects in the relevant derived Fukaya category. As a consequence, at least by standard Floer theoretic methods, it is difficult to distinguish between such $L$ and $L'$.
  \end{rem}

The Proof of Proposition \ref{prop:shadow-est} occupies the rest of this section.

\subsection{Outline of the proof}

We start with an outline of the proof. We first prove the result under the additional assumption that the cobordism $V$ is simple, and that it looks like in Figure \ref{fig:elephant-in-boa}, for a certain function $\beta:\R \to \R$. Later we observe that the proof does not depend on these assumptions. We also assume that $L$
is connected.

Consider a Hamiltonian isotopic copy $L_{\epsilon}\subset M$ of $L$ so that $w(L_{\epsilon}, L')$ is $\epsilon_{0}/2$ close to $w(L,L')$ and $L_{\epsilon}$ is transverse to both $L$ and $L'$.
Choose a symplectic embedding of a ball $e:B_r \to M \setminus L'$ with half-capacity $\frac{\pi r^2}{2} \geq w(L_{\epsilon},L') - \epsilon_0/2$ sending the real part of $B_r$ to $L_{\epsilon}$.  Consider an $\omega$-compatible almost complex structure $J_{-}$ which extends the push-forward of the standard complex structure on $B_r$ by $e$ and let $J$ be an almost complex structure as in the statement of the Proposition.
To prove our statement, it is enough to show that $\mathcal{S}(V) < \delta(V,J)$ implies
$\mathcal{S}(V)\geq w(L,L')-\epsilon_{0},$ which would follow from $\mathcal{S}(V)\geq w(L_{\epsilon}, L')-\epsilon_{0}/2$. We continue to describe the main steps required to prove the last inequality under the assumption  $\mathcal{S}(V) < \delta(V,J)$.  
 
 To show the inequality $\mathcal{S}(V)\geq w(L_{\epsilon},L')-\epsilon_{0}/2$ we will produce a certain $J_{-}$-holomorphic strip $v$ passing through the center of the ball $B_{r}$ and with boundary on $L_{\epsilon}$ and $L'$ whose area is bounded in terms of $\mathcal{S}(V)$.  Given such a strip, usual isoperimetric type inequalities imply that the area of this strip is no less than $\frac{\pi r^{2}}{2}$ and the
 desired inequality follows.
 
The existence of such a $J$-curve will be deduced from a comparison among truncated Floer complexes associated to certain cylindrical Lagrangians constructed out of $V$ as follows. In  (Section \ref{subsubsec:reposition}) we modify $V$ so that its projection is as in Figure \ref{fig:tilde-beta}. Then we perform an additional modification of $V$ to $V'$ as in Figure \ref{fig:bend-V}, so that $\pi(V')$ transversely intersects the $\R+i$-axis in a unique point $P$. The modifications depend on small 
parameters $\delta,\delta',\delta'',\delta''' > 0$, that will be sent to $0$ at the end of the argument. We then construct a special compactly supported, non-negative Hamiltonian $H$ on $\C \times M$. The shadow of $V$ enters the proof through the oscillation of $H$, $\osc(H)$, which is given by $\osc(H)=(1+\delta'')(\mathcal{S}(V)+\delta)$.  The additional key property of $H$ is the following. If we put
 $V''=\phi^{H}_{1}(V')$, then  $\pi(V'')$ intersects the real axis transversely at a unique point $Q$, see   Figure \ref{fig:bend-V}. We also require that $\{\phi_t^H\}_{t \in [0,1]}$ is $i \times J_{-}$-holomorphic in a narrow vertical rectangle $\mathcal{R}$ containing $Q$ and so that the area of $\mathcal{R}$ is  bigger than the oscillation of $H$. Let $L_\epsilon$ denote a Lagrangian submanifold of $M$ that is obtained by a small Hamiltonian perturbation of $L$ and is transverse to both $L$ and $L'$.  Define an additional cylindrical Lagrangian by  \[W = \R \times \{1\} \times L_\epsilon.\] 

We can now be more precise concerning the properties of the curve that we intend to produce:
 it is enough to prove that there exists a $J_{-}$-holomorphic curve $v$ with boundary on $L_{\epsilon}$ and $L'$, that passes through the center $e(0)=R \in L$
of $e(B_{r})$, and that has $\omega$-area at most $\osc(H)+2\zeta$ where the constant $\zeta>0$ can be taken as small as desired. Moreover, this in turn follows if we prove the same statement for a generic almost complex structures instead of $J_{-}$, by virtue of Gromov compactness.

To prove the existence of such a curve $v$ we use a sandwich argument involving 
three Floer complexes (Section \ref{subsubsec:filt-Fl}):
$CF^{a}_{b}(W,V';-H, J)$, $CF^{a}_{b}(W,V';H_{1}, J)$,  $CF^{a}_{b}(W, V'; H_{2}, J)$. Here $H$ is the Hamiltonian above and $H_{1}$ and $H_{2}$ are constant Hamiltonians on $\C\times M$
so that we have $H_{2}\geq - H\geq H_{1}$. The constants $a>b\in \R$ indicate that the respective complexes are taken in the action window $[b,a]$. They are picked so that $ a-b < \delta(V,J)$, $a=\zeta>0$, $b=-\osc(H)-\zeta$.  Further, the constant Hamiltonians are picked so that
 $H_{2}\equiv 0$ and $H_{1}\equiv -\osc(H)-\delta''$ where $0<\delta''<\zeta$.
We use the usual Floer action functional defined on a cover $\widetilde{\mathcal{P}}_0$ of the component $\mathcal{P}_0$ of the space of paths from $W$ to $V'$ that contains the constant path. The cover is associated to the kernel of the morphism $\pi_{1}(\mathcal{P}_{0})\to \R$ given by integrating $\Omega=\omega_{0}\oplus\omega$. 
 Floer homology is taken with coefficients in $\Z_{2}$ and without grading.
 These complexes are related by natural monotone-homotopy comparison morphisms. The fact that we work under
 the bubbling threshhold implies that these truncated complexes are indeed chain complexes and that the comparison morphisms are chain maps. In this overview we neglect the compactness and regularity issues, 
 and assume implicitly that the complexes above are well defined as stated. In the body of the actual proof these points will be addressed.
 
 Suppose that we can use
 the properties of the complexes $CF^{a}_{b}(W,V';H_{i},J)$ to extract the existence of a Floer strip $u$
 associated to the Floer data giving the complex $CF^{a}_{b}(W,V';-H,J)$ so that $u$ passes through the line 
 $\R\times\{1\}\times\{R\}\subset W$ and has energy at most $a-b=\osc(H)+2\zeta$. To this strip $u$ we then associate a new curve $v(s,t) = (\phi^{H}_{t})(u(s,t))$. It is easy to check that this is a Floer strip for the data $(W, V''; 0, J')$ where $J'$ is the  time
 dependent almost complex structure so that $(\phi^{H}_{t})_{\ast} J = J'$. Given that $\pi(W)$ and $\pi(V'')$ only intersect at $Q$ and given the behaviour of $\phi^{H}_{t}$ over the rectangle $\mathcal{R}$ a simple application of the open mapping theorem implies that $v$ is entirely contained in the fiber $\{Q\}\times M$. 
Morevoer, $v$ passes through $R$ and the $\omega$-area of $v$ equals the energy of $u$ and is therefore bounded
 by $\osc(H)+2\zeta$. Thus, this curve $v$ has the required properties, whence proving the statement is reduced to showing the existence of the curve $u$.
 
To produce the curve $u$ as above we proceed as follows (Section \ref{subsec:Fl-hlgy-Morse}).
On all  complexes $CF^a_b(W,V';F,J)$, where $F=H,H_{1},H_{2}$, there is an action (on the homology level) of the Morse complex $CM(f,W)$ of a function $f_0-|\mathrm{Re}(z)|^2$ on $W$.  Here $f_0:L_\epsilon \to \R$ is a Morse function with a unique minimum at $R$. We use on $W$ a metric of the form $ds^{2}+g_{0}$
where $s$ denotes the coordinate along $\R,$ and $g_{0}$ is a metric on $L_{\epsilon}$ such that the pair $(f_{0},g_{0})$ is Morse-Smale.
Moreover, this action homotopy-commutes with the comparison maps.  Denote  by $R'$ the critical point
of $f$ of index $1$ that corresponds to the minimum of $f_{0}$. Notice that the unstable manifold of this
critical point is precisely the line $\R\times \{1\}\times \{R\}$. It follows from the definition of the module action that to prove the existence of the curve $u$ as above, it is sufficient to show that $R' \ast HF^a_b(W,V';H,J) \neq 0$. Let $\psi_{2,1}:HF^{a}_{b}(W,V';H_{2},J)\to HF^{a}_{b}(W,V';H_{1},J)$ be the comparison
map induced by a monotone homotopy between the constant Hamiltonians $H_{2}$ and $H_{1}$. As this
comparison map is natural and as it commutes with the action of $CM(f,W)$ it follows that to
show $R'\ast HF^{a}_{b}(W,V';H,J)\neq 0$ it suffices to show that 
 $R' \ast \mathrm{image}(\psi_{2,1}) \neq 0$. 

We proceed by noticing that, due again to the open mapping theorem and due to the fact that 
$H_{i}$ are constant Hamiltonians,  the complexes $CF^{a}_{b}(W,V'; H_{i}, J)$ are concentrated 
over the fiber $\{P\}\times M$.  Finally, we use a PSS-type construction (Section \ref{subsubsec:PSS}) 
to compare the module action of $CM(f,W)$ on these complexes to the action of the Morse complex 
$CM(f_{0}, L_{\epsilon})$ on a second Morse complex $CM(\sigma, L_{\epsilon})$. The action of 
$R$ on such a complex $CM(\sigma, L_{\epsilon})$ is non-trivial - this reflects the fact that
 through any generic point on $L_{\epsilon}$ there passes one trajectory of $\sigma$ between
the maximum and the minimum of $\sigma$. We deduce from  this that $R' \ast \mathrm{image}(\psi_{2,1}) \neq 0$.

\subsection{Basic setup}

We will first prove the result under the additional assumption that the cobordism 
$V$ is simple (i.e. $k=s=0$) and that $L$ is connected.  We assume this from now on.
Moreover, to fix ideas, we will assume that the constant $\alpha$ in the definition of the negative end of an almost complex structure that is trivial at infinity is $\alpha_{-}=\alpha_{+}=1$. In particular, over $(-\infty, -1]\times \R$ the almost complex structure coincides with $i\times J_{-}$ and over $[1,\infty)\times \R$ it is of the form $i\times J_{z}$, $J_{z}\in \mathcal{J}$.  Further, again to simplify notation we will assume that the given cobordisms are cylindrical outside the region $[0,1]\times \R$, in other words, for the positive ends the constant $\beta_{+}$ is no bigger than $1$ and for the negative ends the constant $\beta_{-}$ is no smaller than $0$ - see \S\ref{subsec:cob}. This condition
is very easy to achieve by using an appropriate horizontal isotopy.

\subsection{Simplifying assumptions and two auxiliary cobordisms.}\label{subsubsec:reposition} 
We start by proving the  statement by assuming that 
$V$ is positioned as in Figure \ref{fig:elephant-in-boa} below. We will see at the end of the proof that 
the same arguments apply to the general case.
\begin{figure}[htbp]
   \begin{center}
      \epsfig{file=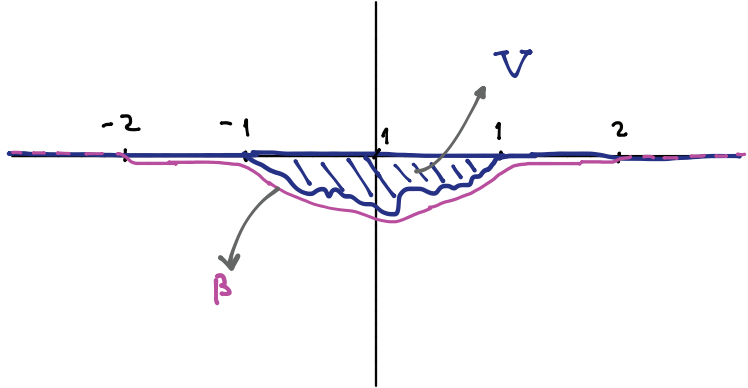, width=0.6\linewidth}
   \end{center}
   \caption{\label{fig:elephant-in-boa} The cobordism $V$ projected onto $\C$ and the graph of $\beta$.
   In this picture as well as the following ones the horizontal axis is at height $1$.}
\end{figure}

\begin{figure}[htbp]
   \begin{center}
      \epsfig{file=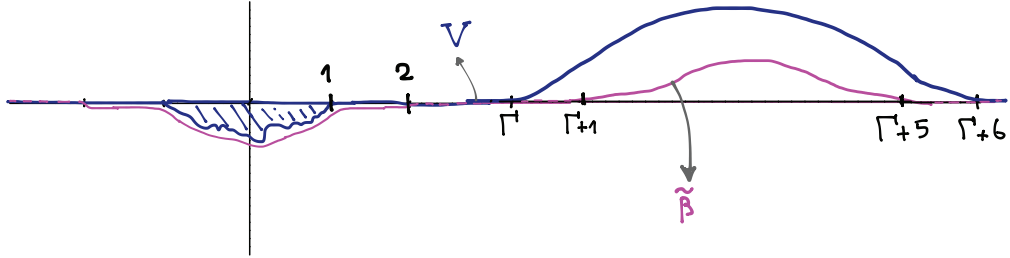, width=0.8\linewidth}
   \end{center}
   \caption{\label{fig:tilde-beta} The graph of $\tilde{\beta}$ and the cobordism $V$ transformed by bending  its positive end.}
\end{figure}

 More explicitly,
the assumption is that:
\begin{itemize} 
\item[i.] $\pi({V}) \cap (\C\setminus [-1,1]\times \R) =  ((-\infty,-1]\cup [1,\infty))\times \{1\}$.
\item[ii.] $\pi({V}) \subset \{ (x+iy)\in \C\ |\ \beta(x) \leq y\leq 1 \ , \ x\in\R\}$ where $\beta: \R\to (-\infty, 1]$
is a smooth function so that the support of $1-\beta(x)$ is inside $[-2,2]$ and $\beta(x)<1$ for $x\in (-2,-1]\cup [1,2)$.
\item[iii.] $\int_{-\infty}^{\infty} (1- \beta (x))dx = \mathcal{S}({V})+\delta$ for some $\delta >0$.
\end{itemize}
 Further, notice that  the constant $\delta$ can be made as small as desired so that
 it is enough to prove:
\begin{equation}\label{eq:mod-est}
\mathcal{S}(V) +\delta \geq \min \{w(L,L')-\epsilon_0, \delta(V,J)\}~,~
\end{equation}
for $J$ as in the statement of Proposition \ref{prop:shadow-est}. Using the function $\beta$ we define a new function $\tilde{\beta}:\R\to \R$ as follows: fix a large 
positive constant $\Gamma >2$ and let $\tilde{\beta}(x)=\beta(x)$ for $x\leq \Gamma$ and 
$\tilde{\beta}(x)=2-\beta(x-\Gamma -3)$ for $x > \Gamma$ - see Figure \ref{fig:tilde-beta}.

Consider now the following function $h_{0}:\C\to \R$:
$$h_{0}(x+iy)=\int_{-\infty}^{x}(1-\tilde{\beta}(t))dt~.~$$
Obviously, this function only depends on $x$ and it is constant equal to $0$ for $x\leq -2$ it is then increasing
till $x=2$ remains constant maximal till $x=\Gamma+1$ when it starts decreasing till it reaches $0$ again for 
$x=\Gamma+5$. It remains equal to $0$ after that. In particular, $h_{0}$ is a positive function and its 
maximum is equal to $\int_{-\infty}^{\infty}(1-\beta(t))dt$.

Let $B>0$ be a constant so that $-B<\min_{x}\beta(x)-2$ and $B\geq \max_{x} (2-\beta(x))+2.$ We remark that $B$ subject to this condition can be chosen to be arbitrarily large - we will use this freedom of choice below.
Consider a new function $h:\C\to  [0,\infty)$  that is supported inside 
$[-4, \Gamma +8]\times [-B-1,B+1]$ and that has the following two properties: $h$ is equal to
$h_{0}$ over the band $\R\times [-B, B]$ and $h\leq h_{0}$ everywhere over $\C$. 
Such an $h$ is easy to construct by cutting off $h_{0}$ outside $\R\times [-B,B]$.

At this point we need to also modify $V$ in a simple fashion that takes into account the ``shape'' of $\tilde{\beta}$.
The modification only consists in ``bending'' the positive end of $V$  by replacing the portion of this end given by
the region $[\Gamma, \Gamma +6]\times \{1\}\times L$ with $l\times L$ where $l$ is the graph of a smooth function $r:[\Gamma, \Gamma +6]\to [1, B)$ with $r(t)>\tilde{\beta}(t),\ \forall \ t$
and $r(t)=1$ for $t\in [\Gamma, \Gamma +\frac{1}{2}]\cup [\Gamma +\frac{11}{2}, \Gamma+6]$.
We will denote the resulting cylindrical Lagrangian still by $V$ - see Figure \ref{fig:tilde-beta}. It obviously has the same shadow as the initial 
cobordism.

We now consider the Hamiltonian
$\bar{h}:\C\times M\to \R$ defined by $\bar{h}=h\circ \pi$ and list its properties that will play an
important role later in the proof:
\begin{itemize}
\item[a.] The oscillation of $\bar{h}$ is  $\int_{-\infty}^{\infty}(1-\beta(x))dx$.
\item[b.] The Hamiltonian vector field $X^{\bar{h}}$ is a horizontal lift of the vector field $X^{h}$
which has compact support.
\item[c.] Over the set $\R\times [-B,B]$,  we have $$X^{h}(x+iy)=i(1-\tilde{\beta}(x))~.~$$
\item[d.] Because of point {\em c} and, due to the perturbation of the positive end of $V$ involving the curve $l$ above, the
 time one diffeomorphism associated to $X^{\bar{h}}$, $\phi^{\bar{h}}$, has the property 
that $\pi(\phi^{\bar{h}}(V))\subset \R\times [1,\infty)$.
\end{itemize} 

\begin{figure}[htbp]
   \begin{center}
      \epsfig{file=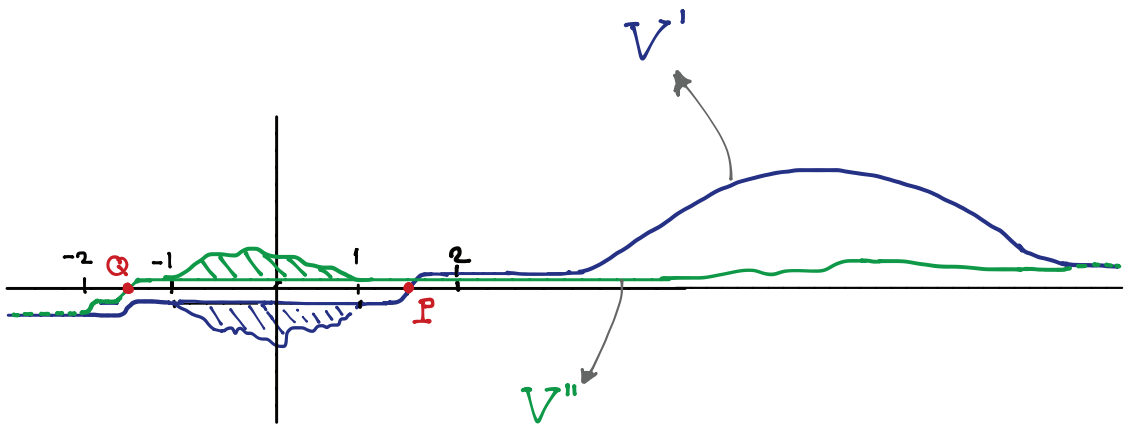, width=0.9\linewidth}
   \end{center}
   \caption{\label{fig:bend-V} The two Lagrangians with cylindrical ends $V'$ and $V''$ projected on
   $\C$.}
\end{figure}

We need to construct two further Lagrangians with cylindrical ends that are basically copies of $V$.
The first will be denoted $V'$. It is obtained by first cutting $V$ along $L_{0}=\{\frac{3}{2}+i\}\times L$
into two pieces one with projection to the left of $(\frac{3}{2}+i)$ denoted by $V_{-},$ and one with projection to the right of $(\frac{3}{2}+i)$, denoted by $V_{+}$.  We translate $V_{-}$ using the 
transformation $(x+iy, m)\to (x+i(y-\delta'), m)$ for $\delta'>0$ very small thus getting a copy of $V_{-}$
denoted by $V'_{-}$. We translate $V_{+}$ using the transformation $(x+iy,m)\to (x+i(y+\delta'),m)$
 and denote by $V'_{+}$ the resulting submanifold with boundary. We then bend slightly the ends at $L_{0}$
 of $V'_{-}$ and $V'_{+}$ and glue them together. We denote by $V'$ the resulting 
 cylindrical  Lagrangian  - as in Figure \ref{fig:bend-V}. 
 We do this construction so that 
 the Lagrangian $V'$  has the property that  $\pi(V')$ intersects $\R\times \{1\}$ transversely, 
in just a single point equal to $P=\frac{3}{2}+i\in \C$. The positive end of $V'$ is
now at height $1+\delta'$ .  The negative end of $V'$ is at height
$1-\delta'$.  Finally, we also need to make another slight adjustement to $V'$. We perturb it so that
that its projection $\pi(V')$ is the graph of an increasing, non-constant function
over the region $[-2,-1]\times \R$ - see Figure \ref{fig:bend-V}.
 
 The second Lagrangian will be denoted by $V''$. To define it take yet another small
constant $\delta''$ and consider the Hamiltonian $H=(1+\delta'')\bar{h}$ together with the associated
 Hamiltonian isotopy $\phi^{H}_{t}$. Now define $$V''=\phi^{H}_{1}(V')~.~$$ 
First notice that $V''$ and $V'$ are horizontally isotopic (in fact, the ends of $V'$ remain
fixed during the isotopy). Further, recall that $\beta(x)<1$ for $x\in (-2,-1]\cup [1,2)$
so that by taking $\delta'$ sufficiently small and possibly adjusting slightly $\beta$ and $V'$ 
we may assume that the only intersection between $\pi(V'')$ and $\R\times \{1\}$ is a single point
$Q\in (-2,-1)\times \{1\}\subset \C$ and that this intersection is transverse. For convenience
we put $q=\mathrm{Re}(Q)\in (-2,-1)$.
There is an additional technical assumption on the function $\beta$ that will be required below.  \begin{itemize}
\item[iv.]There is a small constant $\delta'''>0$ so that the function $\beta$ is constant in the
interval $(q-\delta''',q+\delta''')\subset (-2,-1)$.
\end{itemize}
Due to the behaviour of $\pi(V')$ in the region $[-2,-1]\times \R$ this can be easily achieved by a small perturbation of $\beta$. Finally, we choose the cutoff parameter $B$ in $h$ so that the area $A=2\delta'''\cdot(B-3)$ of each of the two rectangles $\{q-\delta'''<x<q+\delta''',\; 2 < \pm y < B-1\}$ satisfies $A > \osc(H).$

\

To summarize the construction, we have constructed two cylindrical Lagrangians $V'$ and $V''$
that are as in Figure \ref{fig:bend-V}. Moreover, $V'$ and $V''$ both have ends $L$ and $L'$
and they are horizontally Hamiltonian isotopic via the Hamiltonian $H$ whose oscillation is 
$$\osc(H) =(1+\delta'')\int_{-\infty}^{\infty}(1-\beta(t))dt = (1+\delta'')(\mathcal{S}(V)+\delta)~.~$$
To show (\ref{eq:mod-est}) it is enough to prove that
\begin{equation}\label{eq:mod-est2}
\osc(H)\geq \min\{w(L,L')-\epsilon_0, \delta(V,J)\}~.~ 
\end{equation}

\subsection{A comparison Lagrangian and the almost complex structure $J_{-}$.}\label{subsec:comp-Lagr}
For the arguments below we will need a third Lagrangian denoted $W$ that we define now.
We consider a small Hamiltonian deformation $L_{\epsilon}$ of $L$ inside $M$ so that $L_{\epsilon}$ intersects $L$
as well as $L'$ transversely and, moreover,  $w(L_{\epsilon}, L')$ does not differ from $w(L,L')$ by more than $\epsilon_0/2$. 
We may take $L_{\epsilon}$ in a Weinstein neighbourhood of $L$ to be the graph of a form $d\kappa$ where $\kappa:L\to \R$ is a sufficiently small Morse function. With this notation the Lagrangian $W$ is simply 
$$W=\R\times\{1\}\times L_{\epsilon}~.~$$
In Figure \ref{fig:bend-V} the projection of $W$ on the plane $\C$ coincides with the
horizontal line $\R\times \{1\}$. Notice that $W$ intersects transversely both $V'$ and $V''.$ To prove the first point of the proposition for a simple cobordism, it is enough to show:
\begin{equation}\label{eq:mod-est3}
\osc(H)\geq \min\{w(L_{\epsilon},L')-\epsilon_0/2, \delta(V,J)\}~.~ 
\end{equation}
Using the Lagrangians $L_{\epsilon}$ and $L'$ we can now explain our choice of almost complex structure
$J_{-}$. First consider a symplectic embedding of the standard ball $e:B_{r} \to M$ so that $e^{-1}(L_{\epsilon})=B_{r}\cap \R^{n}$, $e(B_{r})\cap L'=\emptyset$ and $$\frac{\pi r^{2}}{2}> w(L_{\epsilon},L')-\epsilon_0/2.$$ Let $R=e(0)$ and let $J_{-}$ be 
an almost complex structure on $M$ that is compatible with $\omega$ and is an extension to the exterior of $e(B_{r})$ of the almost complex structure $e_{\ast} (i)$  ($i$ is here the standard complex structure on $\C^{n}$). 
Claim (\ref{eq:mod-est3}) follows if we can show that, assuming $\osc(H)<\delta(V,J)$, 
\begin{equation}\label{eq:main-cl}
\osc(H)\geq w(L_{\epsilon}, L')-\epsilon_{0}/2~.~
\end{equation}

We pursue the proof under the assumption $\osc(H)< \delta(V,J)$ and with the aim to prove
(\ref{eq:main-cl}).

\subsection{Filtered Floer complexes.}\label{subsubsec:filt-Fl}  Let $J$ be an almost complex structure as in the statemenent of the first point of the 
Proposition \ref{prop:shadow-est}. We will discuss in this subsection the construction of a few truncated Floer complexes together with comparison maps relating them that are the basic tools in the proof.  
These complexes are of the following form:
$$C_{G,J}=CF^{a}_{b}(W,V';G, J)$$
for $a,b\in \R$, $0 < a-b< \delta(V,J)$, and $G:\C\times M\to \R$ a Hamiltonian such that $dG$ is compactly
supported and that $\osc(G)<a-b$ and $J$ is the almost complex structure fixed above.

One point should already be emphasized here: as usual in Floer theory, to insure regularity one might need to
use generic perturbations of the initial data $G,J$ in the construction. For simplicity we do not change the notation at this point but we will discuss this issue in more detail further below.
 
To start the construction of the complex $C_{G,J}$ recall that $L_{\epsilon}$ is the graph of a form $d\kappa$ with $\kappa:L\to \R$ a Morse function. We will assume here that $\kappa$ has a single minimum $m_{0}$ and a single maximum denoted $w_{0}$. Notice that we can view $m_{0}$ as well as $w_{0}$ as (constant) paths from $W$ to $V'$. 
Consider the path space $\mathcal{P}(W, V')$ which consists of the smooth paths $\gamma :[0,1]\to \C\times M$, $\gamma(0)\in W$, 
$\gamma(1)\in V'$. Let
$\mathcal{P}_{0}$ be the path component of $\mathcal{P}(W, V')$ that contains $m_{0}$. Pick $m_0$ as a basepoint in $\mathcal{P}_{0}$. Denote by  $\widetilde{\mathcal{P}_{0}}$ the cover of $\mathcal{P}_{0}$ corresponding to the subgroup $\ker I_\Omega \subset \pi_1(\mathcal{P}_0),$ where $I_\Omega: \pi_1(\mathcal{P}_0) \to \R$ is given by integrating the symplectic form $\Omega = \omega_0 \oplus \omega$. Finally, let $p:\widetilde{\mathcal{P}_{0}}\to \mathcal{P}_{0}$ be the projection. 

As a vector space the complex $CF(W,V';G,J)$ is the $\Z_{2}$ vector space freely generated by the set $\widetilde{\Gamma}_{G}$ of those $\gamma$ in $\widetilde{\mathcal{P}_{0}}$ that project to paths $x = p(\gamma) \in \mathcal{P}_{0}$
that are Hamiltonian chords for $G$ in the sense that they satisfy
$$\frac{d x}{dt}=X^{G}(x)$$
where $X^{G}$ is the Hamiltonian
vector field of $G$.
We denote by $\Gamma_{G}$ all these Hamiltonian chords.
Fix also a lift $\tilde{m}_{0}$ of $m_{0}$ (say as a constant path) to $\widetilde{\mathcal{P}_{0}}$ and define the action
of this point $\tilde{m}_{0}$ by $\mathcal{A}_{G}(\tilde{m}_{0})=G(m_{0})$. Further, for any $\gamma\in \widetilde{\mathcal{P}_{0}}$ define $$\mathcal{A}_{G}(\gamma)=\int_{0}^{1}G((p\circ\gamma)(t))dt -\int_{[0,1]\times [0,1]} (p\circ\bar{\gamma})^{\ast}(\omega)$$
where $\bar{\gamma}$ is a path in $\widetilde{\mathcal{P}}_{0}$
that  joins, in this order, $\tilde{m}_{0}$ to $\gamma$. We will say that $\tilde{m}_{0}$ is the basepoint for the action. 
It is easy to see that the action $\mathcal{A}_{G}(\gamma)$ is independent of the choice of $\bar{\gamma}$. Notice that we work here in the absence of grading and orientations.
The presumtive differential of the Floer complex $d:CF(W,V';G,J)\to CF(W,V';G,J)$ is defined by 
$$d\tilde{x}=\sum_{y}\#_{2}\mathcal{M}(\tilde{x},\tilde{y};G,J)\tilde{y}$$  
where $\tilde{x},\tilde{y}\in \widetilde{\Gamma}_{G}$,
 and $\mathcal{M}(\tilde{x},\tilde{y};G,J)$  is a moduli space consisting
 of paths $\tilde{u}:\R\to \widetilde{\mathcal{P}_{0}}$ - modulo reparametrization by the action of $\R$ -
  that go from $\tilde{x}$ to $\tilde{y}$ and so that the map $u=p(\tilde{u})$ is a Floer strip $u$ from $x=p(\tilde{x})$ to $y=p(\tilde{y})$. In other words $u$
 verifies the Floer equation:
\begin{eqnarray}\label{eq:Floer-strips}
u:\R\times [0,1]\to \C\times M\ ,\ u(\R\times\{0\})\subset W \ , \ u(\R\times \{1\})\subset V' \\
\nonumber \lim_{s\to -\infty}u(s,t)= x(t) \ , \ \lim_{s\to +\infty}u(s,t)=y(t) \ ,\ 
\frac{\partial u}{\partial s} + J(t,u)\frac{\partial u}{\partial t} +\nabla G(u)=0 ~,~
\end{eqnarray}
where the gradient $\nabla G(u) = - J(t,u) X^G(u)$ is computed with respect to the metric $\omega(-,J-)$ given by $J$ and $\omega.$

By definition, $\#_{2}A $ vanishes whenever the set $A$ is infinite and equals the number (mod $2$) of elements
of $A$ when $A$ is finite.  Recall that the almost complex structure is time dependent which explains the notation $J(t,u)$. 
The energy of a strip $u$ as before is defined as 
$$E(u)= \frac{1}{2}\int_{\R\times [0,1]}(||\frac{\partial u}{\partial s}||^{2} + ||\frac{\partial u}{\partial t}-X^{G}(u)||^{2})dsdt$$
where the norms are taken with respect to the metric $\omega(-, J-)$. This energy is equal to the 
difference of actions $$E(u)=\mathcal{A}_{G}(\tilde{x})-\mathcal{A}_{G}(\tilde{y})~.~$$

The square of the linear map $d$ defined as above does not generally vanish. To extract from this definition a useful chain complex we need to address a few issues that are typical for
Floer type constructions:

\begin{itemize}
\item[i.] Compactness. There are three different phenomena that can lead to lack of compactness of the $1$-dimensional moduli spaces:
\begin{itemize}
\item[a.] Breaking along Hamiltonian orbits in $\Gamma_{G}$. This is precisely the type of non-compactness
that appears in showing that $d^{2}=0$.
\item[b.] Bubbling of disks or spheres. We will deal with this issue considering the complex in an appropriate action window, as we explain below following closely \cite{Bar-Cor:Serre}.
\item[c.] The fact that the target manifold $\C\times M$ is not compact.
Outside of a large compact $\widetilde{K}_{G,J}\subset \C$ the gradient $\nabla G$ vanishes and the almost complex structure
$J$ can be written, at each point $z$, as a product $i\times J_{z}$. As a consequence, if $u$ is a solution of (\ref{eq:Floer-strips}) and we put $u'=\pi\circ u$, then, outside of $\widetilde{K}_{G,J}$, the strip $u'$ is holomorphic and 
an easy application of the open mapping theorem (see also \cite{Bi-Co:cob1}) immediately implies that
$u'$ does not get out of $\widetilde{K}_{G,J}$. 
\end{itemize}
\item[ii.]  To achieve the regularity of the moduli spaces considered here we may need to replace the Hamiltonian $G$ and the almost complex structure  $J$ by arbitrarily small perturbations (possibly time dependent). Any generic pertubations $\widetilde{G}$ of $G$ and $\widetilde{J}$ of 
$J$ with the property that  $d\widetilde{G}$ is 
supported inside $\widetilde{K}_{G,J}$ and so that $\widetilde{J}$ continues to have a product form outside of $\widetilde{K}_{G,J}$ allows for the associated Floer complex $C_{\widetilde{G},\widetilde{J}}$ to be defined.
This is the case because our Floer strips $u$ are known {\em apriori} to remain inside $\widetilde{K}_{G,J}\times M$ (by point (i.c.)). Of course, in the constructions below some special choices of such perturbations are sometimes sufficient. They will be discussed explicitly when needed.
\end{itemize}
The truncation mentioned at the point (i.b) above is based on the fact that the action $\mathcal{A}_{G}$ decreases along Floer trajectories. This means that for $a\in \R$ we can define
$$CF^{a}(W,V';G,J)=\Z_{2} <\tilde{x}\in \widetilde{\Gamma}_{G} \ | \ 
\mathcal{A}_{G}(\tilde{x}) < a >$$
and, further, for some other $b\in \R$, $b\leq a$,
$$CF^{a}_{b}(W,V';G,J)=CF^{a}(W,V';G,J)/CF^{b}(W,V';G,J)~.~$$
This definition is of interest as soon as $a-b < \delta(V,J)$ as this implies that the linear map induced by $d$, as defined through (\ref{eq:Floer-strips}), is a differential on $CF^{a}_{b}(W,V';G,J)$. To see this we only need to notice that for appropriate choices of parameters in the construction of $V'$ and $W,$ any potential $J$-sphere or disk that can bubble off in a $1$-dimensional moduli space has symplectic area at least $\delta(V,J)-\xi > a-b$ for some very small positive constant $\xi$. 
At this point it is useful to recall that $J$ is time dependent and to focus on what type of such curves can potentially bubble off: these are $J_{t}$-spheres for any $t\in [0,1]$ - by the definition of $\delta(V,J)$ they have an area no less than this number; $J_{1}$-disks with boundary on $V'$ - by Gromov compactness, these are of area at least $\delta(V,J)-\xi$ as soon as $V'$ is sufficiently close to $V$ and if we take the constant  $\Gamma$ in the ``bending'' perturbation of $V$ big enough so that the bend takes place in the region where $J$ is of the form $i\times J_{+}$;  finally, there can also be $J_{0}$-disks with boundary on $W$.  Recall that $W=\R\times\{1\}\times L_{\epsilon}$ and $J_{0}=i\times J_{-}$ so that that the symplectic area of these $J_{0}$ disks is at least $\delta(L_{\epsilon}, J_{-})$. By Gromov compactness, taking $L_{\epsilon}$ close enough to 
$L$, we have $\delta(L_{\epsilon}, J_{-})\geq \delta(L,J_{-})-\xi$.  But $L$ is also the negative end
of $V$ and $J_{-}$ is the negative end of $J_{t}$ for all $t$. It follows that $\delta(L,J_{-})\geq \delta(V,J),$ concluding this argument.  In case $J$ is required to be perturbed for regularity purposes, the resulting
disks and spheres will still have areas at least $\delta(V,J)-\xi$ if the perturbation is small enough. 

For the continuation of the argument we fix  the constant $\xi$  above and we put  
$$\delta'(V,J)=\delta(V,J)-\xi < \osc(H)$$ where $H$ is the particular Hamiltonian constructed in \S\ref{subsubsec:reposition}.
Fix $\zeta>0$ so that 
\begin{equation}\label{eq:est-z}
\osc(H)+2\zeta < \delta'(V,J) ~.~
\end{equation}

With these conventions we now can fix the constants  $a$ and $b$ that we use further in the proof.
 \begin{equation}\label{eq:bounds}
 a=\zeta  \ , \ b=-\osc(H)-\zeta ~.~
 \end{equation} 

Finally, we can define the three truncated Floer complexes that are at the center of the argument.
They are $$C=CF^{a}_{b}(W,V';-H,\widetilde{J}) \ \ \mathrm{and}\  \ C_{i}=CF^{a}_{b}(W,V'; H_{i},\widetilde{J})$$ for $i=1,2$ where the Hamiltonians $H_{i}$ are as follows: $H_{1}$ is the constant Hamiltonian $$H_{1}\equiv -\osc(H)-\delta'' \ \mathrm{with} \ \  \zeta>\delta''>0$$ and $H_{2}$ is the constant Hamiltonian equal to $0$.
The almost complex structure $\widetilde{J}$ is a small generic perturbation of $J$ as we now explain. 
First, we slightly perturb $J_{-}$ to a possibly time dependent almost complex structure $\widetilde{J}_{-}$ on $M$ and $J_{+}$ (which is, in general, already time-dependent) to a 
structure $\widetilde{J}_{+}$ so that the data $(L_{\epsilon}, L; 0 , \widetilde{J}_{+})$ and 
$(L_{\epsilon}, L'; 0,\widetilde{J}_{-})$ are Floer regular. We then extend this perturbation of the ends of $J$ to a generic perturbation $\widetilde{J}$ of $J$ itself that remains in the class of almost complex structures trivial at infinity and so that $\widetilde{J}$ has $\widetilde{J}_{-}$ as negative end  and has $\widetilde{J}_{+}$ as positive end and the associated constants for these ends remain 
$\alpha_{-}=1=\alpha_{+}$.
 Given that the intersection between $W$ and $V'$ is transverse this is enough for the 
regularity required to define the complexes $C_{i}$. Similarly, because the intersection of $W$ and $V''$ is transverse this type of perturbation is also sufficient to define $C$.

\

\subsection{Reduction of (\ref{eq:main-cl}) to the existence of certain Floer strips.}\label{subsec:reduct}
The purpose of this subsection is to notice our claim (\ref{eq:main-cl}) is implied by the following
statement: 
\begin{eqnarray} \label{eq:mod-est6b}\exists \ u:\R\times [0,1]\to \C\times M, \  u  \ \mathrm{verifies\ (\ref{eq:Floer-strips})} \ \mathrm{with} \  G=-H, J=\widetilde{J} \\ \mathrm{and}\  u(0,0)\in \R\times \{1\}\times R\ , \
\nonumber E(u)\leq a-b = \osc(H)+2\zeta\hspace{0.3in}
\end{eqnarray}

Indeed, given such a $u$ let $v(s,t)=\phi^{H}_{t}(u(s,t))$. It is easy to see that
$v$ is a solution of 
\begin{eqnarray}
\label{eq:mod-est5} \ v:\R\times [0,1]\to \C\times M,\ v(\R\times \{0\})\subset W,\  v(\R\times \{1\})\subset V'', \  \bar{\partial}_{J'}v=0 \ \mathrm{and}\\  v(0,0)\in \R\times \{1\}\times R\ , \
\nonumber\Omega(v)\leq \osc(H)+2\zeta \hspace{0.8in}
\end{eqnarray}
for an almost complex structure $J'=(\phi_{t}^{H})_{\ast}(\widetilde{J})$. 

Due to condition (iv) in the construction of the Hamiltonian $H$  in
\S \ref{subsubsec:reposition}, the almost complex structure $J'_{t}$ is of the form $i\times \widetilde{J}_{-}$
in a rectangle of the form $(q-\delta''', q+\delta''')\times [-B+1,B-1]$ where $Q=(q,1)$ is the only point of 
intersection of $\pi(V'')$ and $\pi(W)$. This is because over that strip $\phi^{H}_{t}$ is just a translation
in the imaginary direction in this rectangle (in $\C$). Consider now a curve $v$ given as in (\ref{eq:mod-est5})
and let $v'=\pi\circ v$. This is a curve that is holomorphic over  $(q-\delta''', q+\delta''')\times [-B+1,B-1]$, exits $Q$ and also ends at $Q$. By first using an orientation argument we deduce that this
is not possible if $v'$ does not pass through one of the regions of $[(q-\delta''', q+\delta''')\times \R]\setminus [\pi(W)\cup \pi(V'')]$. But due to the open mapping theorem, and the condition $2\delta'''\cdot (B-3) > \osc(H),$ this is again not possible (for reasons of symplectic area). We deduce that $v'$ is constant and equal to $Q$ and therefore $v$ is contained in $\{Q\}\times M,$ and is a solution of:
\begin{eqnarray}
\label{eq:mod-est4} \ w:\R\times [0,1]\to M,\ w(\R\times \{0\})\subset L_{\epsilon},\  w(\R\times \{1\})\subset L', \  \bar{\partial}_{\widetilde{J}_{-}}w=0 \ \mathrm{and}\\  R=w(0,0)\ , \
\nonumber\omega(w) = \int_w \omega\leq \osc(H) +2\zeta.\hspace{0.8in}
\end{eqnarray}
The  almost complex structure $\widetilde{J}_{-}$ can be taken as close as needed to $J_{-}$ so that, 
by Gromov compactness, we deduce the 
existence of a holomorphic curve $w'$ passing through the center of $B_{r}$ with boundary
on $\partial B_{r}$ and on $\R^{n}\cap B_{r}$ and with $\omega(w')\leq \omega(w)$. The
existence of $w'$ implies, by a standard argument based on the Lelong inequality (see \cite{GriffithsHarris}),
the inequality $$\omega(w')\geq \frac{\pi r^{2}}{2}~.~$$ 
As $\zeta$ can be taken as close to $0$ as desired this implies the  inequality (\ref{eq:main-cl}).

\

Thus, from now on our aim is to show (\ref{eq:mod-est6b}).

\subsection{Comparison of the complexes $C, C_{1}, C_{2}$.}\label{subsec:compare}
Returning now to our Hamiltonians,  notice that $H_{2} \geq -H \geq H_{1}$.
As usual in Floer theory, to compare two Floer complexes associated to different Hamiltonians, we use homotopies that relate these Hamiltonians. In our case, because we deal with truncated complexes it is important to use (decreasing) monotone homotopies: these
 ``reduce'' the action and thus preserve the filtration. 
 There are monotone homotopies
 $G_{2}:[0,1]\times \C\times M\to \R$ from $H_{2}$ to $-H$ and 
$G_{1}:[0,1]\times \C\times M\to \R$ from $-H$ to $H_{1}$ that have the property that for each value of 
the parameter $s\in [0,1]$, the function $G_{i}(s,-)$ is constant outside the large compact $\widetilde{K}_{H,J}\subset \C$ (as described at the point i.c above). Monotonicity means that $\partial_{s} G_{i}\leq 0$.
This means that the linear maps 
\begin{equation}\label{eq:comp1}\psi_{G_{2}}: C_{2}=CF^{a}_{b}(W,V'; H_{2},\widetilde{J})\to CF^{a}_{b}(W,V';-H, \widetilde{J})=C\end{equation}
and
\begin{equation}\label{eq:comp2}\psi_{G_{1}}: C=CF^{a}_{b}(W,V'; -H,\widetilde{J})\to CF^{a}_{b}(W,V';H_{1}, \widetilde{J})=C_{1}
\end{equation}
 defined by counting solutions of an equation similar to (\ref{eq:Floer-strips}): $\bar{\partial}_{\widetilde{J}}u(s,t) +\nabla G_{i}(s,u(s,t))=0$
(for convenience we assume here that $G_{i}$ has been extended to $\R\times \C\times M$ by $G_i(0,-)$ for $s<0$ and  by $G_i(1,-)$ for $s>1$) are chain maps.  Indeed, the same reasons that insure compactness for the solutions of (\ref{eq:Floer-strips}) apply here too. In what concerns regularity we can again perturb in a time 
dependent way the homotopies $G_{i}$ as well as use an appropriate generic $s$-dependent homotopy from $\widetilde{J}$ to $\widetilde{J}$.  Chain maps defined in this way are called  comparison maps in Floer theory and they have an additional naturality property that is important for us:  the
composition $\psi_{G_{1}}\circ \psi_{G_{2}}$ is chain homotopic to the comparison map $\psi_{G_{1}\circ G_{2}}$
and, in turn, this chain map is chain homotopic to any comparison map induced by a monotone homotopy 
from $H_{2}$ to $H_{1}$. This fact is again a consequence of the general machinery involving monotone 
homotopies (see for instance \cite{Bar-Cor:Serre}) together with the same type of arguments insuring compactness and regularity as discussed above.
We denote by $$\psi_{2,1}: CF^{a}_{b}(W,V';H_{2},\widetilde{J})\to CF^{a}_{b}(W,V';H_{1},\widetilde{J})$$ this comparison map.

\subsection{Floer homology as a module over Morse homology}\label{subsec:Fl-hlgy-Morse}   For two Lagrangian submanifolds $N,N'$, there is a module multiplication making $HF(N,N')$ a module over the Morse homology of $N$  viewed as an algebra with the intersection product. Structures of this type appear often in the literature,
 see in particular \cite{Cha:thesis} where the author uses this action for a purpose similar to our aim
here - to detect Floer trajectories through a point. In our application we will actually only need to use
some very basic properties of the module product.  Another technique to detect Floer trajectories through
a point appeared in \cite{Bar-Cor:Serre} but the module multiplication approach is simpler to implement 
here.
 
 Let  $f:W\to \R$ be a Morse function constructed as follows. Recall that $W=\R\times \{1\}\times L_{\epsilon}$. Let $f_{0}:L_{\epsilon}\to \R$ be a fixed Morse function with a single maximum and a single minimum and put $f(s,l)=-s^{2}+ f_{0}(l)$, $s\in\R$, $l\in L_{\epsilon}$.   
 Fix also a Riemannian metric $g_{0}$ on $L_{\epsilon}$ so that the pair $(f_{0},g_{0})$ is Morse-Smale. 
 Extend the metric $g_{0}$ to the metric $g= ds^{2}\oplus g_{0}$ on $W$ and denote by $\gamma_{t}$ the negative
 gradient flow of $f$ with respect to $g$.
 
  Denote by $CM(f,W)$
 the Morse complex associated to $(f,g)$. Note that the obvious map $CM(f_0,L_{\epsilon}) \to CM(f,W)$ is an isomorphism (ignoring grading). Consider also the complex $C=CF^{a}_{b}(W,V';-H,\widetilde{J})$. 
 Notice that there is a module multiplication map:
 $$\mu_{C}:CM(f,W)\otimes C \to C\ , \ \mu_{C}(a,\tilde{x})=\sum_{y}\#_{2}\mathcal{M}(a;\tilde{x},\tilde{y};-H,\widetilde{J})\tilde{y}$$ where $a\in Crit(f)$, $\tilde{x}\in \widetilde{\Gamma}_{H}$ and 
  $\mathcal{M}(a;\tilde{x},\tilde{y};-H,\widetilde{J})$ is the moduli space of 
 paths $\tilde{u}:\R\to \widetilde{\mathcal{P}_{0}}$ that verify (\ref{eq:Floer-strips}) and that additionally satisfy the relation $\lim_{t\to-\infty} p(\gamma_t(\tilde{u}(0,0)))=a$.  By the same methods that were discussed in the 
 last section it follows without difficulty that the map $\mu_{C}$ is a chain map (in particular working under the bubbling threshold is important here).
 
Similarly, we have corresponding 
 chain maps $\mu_{C_{i}}$, $i=1,2$ associated to the complexes $C_{i}$
 that are defined in similar ways as above. These maps are related in the obvious sense
 through the comparison maps $\psi_{G_{i}}$. To be more explicit, we claim that
 $\mu_{C}(id_{CM}\otimes \psi_{G_{2}})\simeq \psi_{G_{2}}(\mu_{C_{2}})$
 and similar identities up to chain homotopy for the other comparison map (see \cite{Bi-Co:qrel-long,Bi-Co:cob1} for related analysis). 
 
 \
 
 Using these module multiplications we can reduce our claim  (\ref{eq:mod-est6b}) to
 an algebraic identity, as follows.  We first pick the Morse
function $f_{0}:L_{\epsilon}\to \R$ above so that its minimum point is $R\in L_{\epsilon}$ that appears in (\ref{eq:mod-est6b}).
Notice that, by the definition of the function $f:W\to \R$ and due to the definition of the metric $g$ on $W$,  the unstable manifold of the critical point $R'=\{0\}\times\{1\}\times \{R\}\in W$  of $f$ is precisely the line $\R\times\{1\}\times \{R\}$. 
In view of these choices, it is immediate to see that (\ref{eq:mod-est6b}) is implied by the following non-vanishing of the module multiplication 
 \begin{equation}\label{eq:mod-est6}
 R'\ast HF^{a}_{b}(W,V';-H,\widetilde{J})\not=0
 \end{equation}
 
Further, by using the comparison maps $\psi_{G_{i}}$ and their compatibility with the multiplications
$\mu_{C},\mu_{C'_{i}}$ we see that (\ref{eq:mod-est6}) is implied by:
\begin{equation}\label{eq:mod-est7}
 R'\ast \mathrm{image}(\psi_{2,1}) \not=0
 \end{equation}
 where we recall that $\psi_{2,1}: C_{2}\to C_{1}$ is the monotone comparison map described in
 \S\ref{subsec:compare}, and we use the same notation for the map it induces on homology.

\subsection{Reduction to an identity in the fibre over $P$}\label{subsubsec:red-to-fibre} 
We start by making more precise the choice of the function $\kappa:L\to \R$  from \S\ref{subsec:comp-Lagr}. Recall that this function has the property that the graph of $d\kappa$ is the Lagrangian $L_{\epsilon}$.  Recall also that the point $\tilde{m}_{0}$ is the base point for the actions we use here and that $m_{0}$ is the minimum point of  $\kappa$ and $w_{0}$ is the maximum point of $\kappa$. We assume that $\kappa$ is small enough, and that 
\begin{equation}\label{eq:bounds-q}
\kappa(m_{0})=0  \  \mathrm{and}\  \kappa(w_{0}) < \zeta-\delta'' ~.~
\end{equation} 
 
 Notice that picking $\kappa$ in this way is possible because the choice of the
 constants $\zeta$, $\delta''$ is independent of the choice of $L_{\epsilon}$.
 
 \
 
 We will consider truncated Floer complexes of 
 the form $CF^{a}_{b}(L_{\epsilon}, L;\eta, \widetilde{J}_{+})$. Here $\eta$ is a constant $\eta\in \R$
 and in our argument it will only take two values:   $\eta_{1}=-\osc(H)-\delta''$ and $\eta_{2}=0$.
 
 The construction of this truncated Floer complex follows the same procedure as in \S\ref{subsubsec:filt-Fl} but it is useful to identify here the precise path space in use. 
 Consider the component $\mathcal{P}_0(L_\epsilon,L)$ of $m_0$ in the space $\mathcal{P}(L_\epsilon,L)$ of paths (in $M$) from $L_\epsilon$ to $L,$ and the  inclusion $j:\mathcal{P}_0(L_\epsilon,L) \to \mathcal{P}_{0}$ induced by $\{P\} \times M \hookrightarrow \C \times M.$ Recall that $P$ is the only point of intersection between $\pi(W)$ and $\pi(V')$.
 Let $\widetilde{\mathcal{P}}_{0}(L,L_{\epsilon})$ be the pull-back of the covering 
 space $\widetilde{\mathcal{P}}_{0}\to \mathcal{P}_{0}$ (recall
 that $\mathcal{P}_{0}$  is the component of $m_{0}$ of the of space of paths joining $W$ to $V'$
 and $\widetilde{\mathcal{P}}_{0}$ is the covering of $\mathcal{P}_{0}$ associated to the morphism
 induced by integration of $\Omega$). We use the path space $\widetilde{\mathcal{P}}_{0}(L,L_{\epsilon})$ to construct $CF^{a}_{b}(L_{\epsilon}, L;\eta, \widetilde{J}_{+})$ the point $\tilde{m}_{0}$, which belongs to $\widetilde{\mathcal{P}}_{0}(L,L_{\epsilon})$, taken as the base point for the relevant action. 
We emphasize that our construction implies that $\mathcal{A}_{\eta}(\tilde{m}_{0})=\kappa(m_{0})+\eta=\eta$. 

Given that over $[1,\infty)\times \R$ the almost complex structure $\widetilde{J}$ is of the form $i\times \widetilde{J}_{+}$, an immediate application of the open mapping theorem (as in \cite{Bi-Co:cob1}) implies that 
  \begin{equation}\label{eq:mod-est8}
 \nonumber CF^{a}_{b}(L_{\epsilon}, L;\eta, \widetilde{J}_{+})=CF^{a}_{b}(W,V';\eta,\widetilde{J})~,~
 \end{equation}  
 where the complex $CF^{a}_{b}(W,V';\eta,\widetilde{J})$ is the complex constructed in \S\ref{subsubsec:filt-Fl} for $G$ the constant Hamiltonian $G\equiv\eta$ defined on $\C\times M$.

 It is easy to see that this identification is again compatible with the multiplications involved  so that (\ref{eq:mod-est7}) becomes:
 \begin{equation}\label{eq:mod-est8}
 R\ast \mathrm{image}(\psi_{\eta_2,\eta_1})\not=0,
  \end{equation}
  where $$\psi_{\eta_{2},\eta_{1}}: CF^{a}_{b}(L_{\epsilon}, L; \eta_{2},\widetilde{J}_{+})\to 
  CF^{a}_{b}(L_{\epsilon}, L;\eta_{1},\widetilde{J}_{+})$$ is the 
  comparison morphism associated to a monotone homotopy relating the two constant hamiltonians
  $\eta_{2}=0$ and $\eta_{1}= -\osc(H)-\delta''$ defined on $M$ (we use the same notation for the map on the homology level).

\
 
 We can further simplify this equation by taking into account that $L$ and $L_{\epsilon}$ are Hamiltonian isotopic. Indeed, recall from \S\ref{subsec:comp-Lagr}, \S\ref{subsubsec:red-to-fibre} that $L_{\epsilon}$ is the graph of the form $d\kappa$. In particular,
 there is a Hamiltonian $\bar{\kappa}: M\to \R$ of oscillation equal to the oscillation of the function $\kappa$ so that $\phi^{\bar{\kappa}}_{1}(L)=L_{\epsilon}$. Moreover, on a Weinstein neighobourhood of $L$, $\bar{\kappa}$ has the form $\bar{\kappa}=\kappa\circ p_{L}$
 where $p_{L}$ is the projection on the base on that neighbourhood. We shall use this Hamiltonian and a naturality type transformation to transform equation (\ref{eq:mod-est8}) into an equation only involving the Floer theory of the pair $(L_{\epsilon}, L_{\epsilon})$. The only subtelty is that we also need to transform
 appropriatedly the covering $p':\widetilde{\mathcal{P}}_{0}(L_{\epsilon}, L)\to \mathcal{P}_{0}(L_{\epsilon}, L)$ into a covering $$p'':\widetilde{\mathcal{P}}_{0}(L_{\epsilon}, L_{\epsilon})\to \mathcal{P}_{0}(L_{\epsilon}, L_{\epsilon})~.~$$
 For this purpose consider the Hamiltonian  $\overline{K}$ on $\C \times M$ given as $0 \oplus \overline{\kappa}$ cut off to $0$ away from a neighbourhood of $\{P\} \times M.$ 
 Let $\Psi'$ be the transformation $$\Psi': \gamma(t) \to \phi^{\overline{K}}_{t}(\gamma(t))~.~$$
 
 Put $a_{0}=\Psi'(m_{0})$ and notice that this is actually a constant path (because $m_{0}\in Crit(\kappa)$). 
 Denote by $V'''=\phi^{\overline{K}}_{1}(V')$ and let $\mathcal{P}_{0}(W,V''')$ be the component
 of the space of path from $W$ to $V'''$ that contains $a_{0}$. Let, as usual, 
 $\tilde{p}:\widetilde{\mathcal{P}}_{0}(W,V''')\to \mathcal{P}_{0}(W,V''')$ be the covering space associated to the kernel of the morphism $\pi_{1}(\mathcal{P}_{0}(W,V'''))\to \R$ given by integrating $\Omega$. Finally, let $\mathcal{P}_{0}(L_{\epsilon}, L_{\epsilon})$ be the component of the space of paths from 
 $L_{\epsilon}$ to itself (in $M$) that contains $a_{0}$. 
 The covering $p''$ is the pullback of the covering $\tilde{p}$ by the inclusion
 $\mathcal{P}_{0}(L_{\epsilon},L_{\epsilon})\to  \mathcal{P}_{0}(W,V''')$ induced by $\{P\}\times M\subset \C\times M$. The transformation $\Psi''$ defines a homeomorphism that relates the 
 two coverings $p'$ and $p''$. Denote by $\tilde{a}_{0}$ the image of $\tilde{w}_{0}$.
 
 Now define the Floer complex $CF^{a}_{b}(L_{\epsilon},L_{\epsilon}; \bar{h}, \bar{J})$ by
 the same procedure as in \S\ref{subsubsec:filt-Fl} but by using as base-point $\tilde{a}_{0}$ and
 the path spaces given by the covering $p''$. Here $\bar{h}$ is a Hamiltonian on $M$. 
 The key remark at this point is that the map $\Psi'$ induces an identification:
 
 \begin{equation}\label{eq:nat2}
\Psi':CF^{a}_{b}(L_{\epsilon}, L;\eta, \widetilde{J}_{+})\to CF^{a}_{b}(L_{\epsilon}, L_{\epsilon};\eta + \bar{\kappa},\bar{J})\end{equation}
 where $\bar{J}=(\phi^{\bar{\kappa}}_{t})_{\ast}(\widetilde{J}_{+})$ and $\eta\in \R$.
 
 This map is also compatible with the action of $CM(f_{0}, L_{\epsilon})$ so that (\ref{eq:mod-est8}) becomes:
 \begin{equation}\label{eq:mod-est10}
   R\ast \mathrm{image}(\bar{\psi}_{2,1})\not=0 \end{equation}
 where $$\bar{\psi}_{2,1}:CF^{a}_{b}(L_{\epsilon}, L_{\epsilon};\eta_{2}+\bar{\kappa},\bar{J})\to
 CF^{a}_{b}(L_{\epsilon}, L_{\epsilon};\eta_{1}+\bar{\kappa},\bar{J})$$ is the natural comparison
 map and, we recall, $\eta_{2}=0$ and $\eta_{1}=-\osc(H)-\delta''$.
 
 We have reduced our argument to showing (\ref{eq:mod-est10}) which will be done in the next
 subsection.
  
\begin{rem} It is tempting to directly argue that
the complex  $CF^{a}_{b}(L_{\epsilon},L_{\epsilon}; \eta+\bar{\kappa}, \bar{J})$ 
can be identified with a Morse type complex by reasoning like in  Floer's work \cite{Fl:Witten-complex}. 
However, we are not here in an exact setting and all such identifications
are very sensitive to the action window. In the next subsection we will compare the complex
$CF^{a}_{b}(L_{\epsilon},L_{\epsilon})$ with an appropriate Morse complex but the argument is more involved.
\end{rem}

 \subsection{A PSS type comparison argument and the proof of (\ref{eq:mod-est10})}\label{subsubsec:PSS}
 We start with the remark that for any Morse function $\sigma:L_{\epsilon}\to \R$ so that the pair $(\sigma,g_{0})$ is Morse-Smale, and the pairs $(\sigma,g_0)$ and $(f_0,g_0)$ are suitably in general position with respect to one another, there is a multiplication $CM(f_{0},L_{\epsilon})\otimes CM(\sigma,L_{\epsilon})\to CM(\sigma,L_{\epsilon})$
 so that $R\ast HM(\sigma,L_{\epsilon})$ is non-trivial. This is because its induced product in homology is identified with the singular intersection product which has a unit (given by the fundamental class).
 The next step is to compare $CM(\sigma,L_{\epsilon})$ and $CF^{a}_{b}(L_{\epsilon},L_{\epsilon};\eta + \overline{\kappa}, \bar{J})$ by
 means of the so-called PSS \cite{PSS} maps.  Such comparisons have often been used in the literature before, for instance in \cite{Alb:PSS}\cite{Bar-Cor:NATO} . 
 
 It is necessary to work with a certain extension of the complex $CM(\sigma,L_{\epsilon})$ that takes into
 account the path spaces used to define the Floer complex $CF(L_\epsilon, L_{\epsilon})$. 
 To define this extension, let $p''':\widetilde{L}_{\epsilon}\to L_{\epsilon}$ be the pull-back covering induced from $p''$ by the map $j_{L_{\epsilon}}:L_{\epsilon}\to \mathcal{P}_{0}(L_{\epsilon},L_{\epsilon})$ which sends each point in $L_{\epsilon}$ to the constant path. Denote by
 $CM(\sigma,\widetilde{L}_{\epsilon})$ the  obvious lift of the Morse complex of $\sigma$.
 In other words, this is the Morse complex of $\sigma\circ p'''$: the generators are the lifts
 of the critical points of $\sigma$ and the connecting trajectories are paths in $\widetilde{L}_{\epsilon}$ that project to negative gradient trajectories of $\sigma$. 
 
\

 It is useful to undestand the complex  $CM(\sigma,\widetilde{L}_{\epsilon})$ better. First of all notice that $\tilde{a}_{0}\in \widetilde{L_{\epsilon}}$. Each point in $\widetilde{L_{\epsilon}}$ is identified with a pair formed by a point $z$ in $L_{\epsilon}$ together with a ``weight" $\Omega(z)$ given by the integral of $\Omega$ over a path in $\mathcal{P}_{0}(W,V''')$ that starts at $a_{0}$ and ends at $z$. Obviously, the integral of $\Omega$ along
 any path in $\mathcal{P}_{0}(W,V''')$ that is completely in $L_{\epsilon}$ vanishes. Therefore, the Morse differential does not modify this weight. As a consequence, if we only look at the generators of $CM(\sigma,\widetilde{L}_{\epsilon})$ that are of weight $0$, they form a subcomplex $CM(\sigma,\widetilde{L}_{\epsilon};0)$ which is actually a factor of $CM(\sigma,\widetilde{L}_{\epsilon}),$ and is obviously isomorphic to $CM(\sigma,L_\epsilon).$ Furthermore, the product $R\ast HM(\sigma, \widetilde{L}_{\epsilon};0)$ obviously continues to be non-trivial. 
 
 In summary, to show our claim (\ref{eq:mod-est10}) it is enough to construct, for $\eta\in [\eta_{1}, \eta_2],$ two maps
  $$\phi_{\eta}  :CM(\sigma, L_\epsilon) = CM(\sigma,\widetilde{L}_{\epsilon},0)\to CF^{a}_{b}(L_{\epsilon},L_{\epsilon};\eta+\bar{\kappa}, \bar{J})$$
 and 
 $$\phi'_{\eta}  :CF^{a}_{b}(L_{\epsilon},L_{\epsilon};\eta +\bar{\kappa}, \bar{J})\to CM(\sigma, \widetilde{L}_{\epsilon};0)$$
 show that they are compatible with the multiplication with $CM(f_{0},L_{\epsilon})$,
 are such that $\phi'_{\eta} \circ \phi_{\eta}$ is chain homotopic to the identity and moreover
 $\bar{\psi}_{2,1}\circ \phi_{\eta_{2}}$ is chain homotopic to $\phi_{\eta_{1}}$.
 
 \
 
To simplify notation we shall denote the Hamiltonian $\eta +\bar{\kappa}$ by $F$ and we assume $\eta\in [\eta_{1},\eta_{2}]$.  We also denote the projections in the respective covering spaces by $p$. The construction of $\phi$ is based on counting trajectories $(u,\gamma)$ where $u:\R\to \widetilde{\mathcal{P}}_{0}(L_{\epsilon},L_{\epsilon})$,
$\gamma:(-\infty,0]\to \widetilde{L}_{\epsilon}$ and if we put $u'=p(u)$, $\gamma'=p(\gamma)$,
then we have:
$$u'(\R\times\{0,1\})\subset L\ , \ \partial_{s}(u')+ \bar{J}_{P}(u')\partial_{t}(u')+\theta(s)\nabla
F(u')=0\ , u(+\infty)=y$$ and
$$
\frac{d\gamma'}{dt}=-\nabla \sigma(\gamma')\ ,\  \gamma(-\infty)=x \ , \
\gamma(0)=u(-\infty)~,~$$
 where $x$ is a lift of a critical point of $\sigma$, with $\Omega(x)=0$, and $y$ is a generator of $CF^{a}_{b}(L_{\epsilon},L_{\epsilon};\eta +\bar{\kappa}, \bar{J})$; $\theta$ is a smooth
 cut-off function which is increasing and vanishes for $s\leq 1/2$ and equals
 $1$ for $s\geq 1$. 

The energy of such an element $(u,\gamma)$ is defined in the obvious way by
$E(u,\gamma)=\int ||\partial_{s}u'||^{2}ds dt$ and it is easy to see
that:
$$E(u,\gamma)=I(u) +
\int_{\R\times [0,1]} (u')^{\ast}\omega-\int_{0}^{1}F(y(t))dt$$
where $I(u)=\int_{\R\times [0,1]} \beta'(s)F(u'(s))dsdt$.
The energy verifies
$$E(u,\gamma)=I(u) -\mathcal{A}_{F}(y)\leq \sup (F) -\mathcal{A}_{F}(y)~.~$$
Recall now from (\ref{eq:bounds}) that $a=\zeta$ and $b=-\osc(H)-\zeta$
and we also have from (\ref{eq:bounds-q}):
$$\sup(F)=\eta + \sup (\bar{\kappa}) < \eta + \zeta -\delta'' \leq \zeta =a~.~$$

Therefore, $\mathcal{A}_{F}(y)< a$ so that $\phi$ is well defined. We also need to notice that our definition
keeps the energy under the bubbling threshold so that this map is a chain map. As $CF^{a}_{b}=CF^{a}/CF^{b}$
it follows that the only orbits of interest have action in between $[-\osc(H)-\zeta,\zeta]$, therefore by (\ref{eq:est-z})
we have
$E(u,\gamma)\leq \osc(H)+2\zeta < \delta'(V,J)$.

The construction of the map $\phi'$ is similar. 
We consider orbits that join lifts of Hamiltonian orbits to lifts of critical points  of $\sigma$, again of weight $0$,
except that the pairs $(u,\gamma)$ considered here, start as semi-tubes and end as flow lines of $\sigma$.
The equation verified by $u$ is similar to the one before but instead of the cut-off function
$\theta$ we use the cut-off function $1-\theta$. The energy estimate in this
case gives $$E(u,\gamma)\leq \mathcal{A}_{F}(y)-\inf (F)~.~$$
Thus $E(u,\gamma)\leq \mathcal{A}_{F}(y)-\eta\leq \osc(H)+\zeta +\zeta < \delta'(V,J)$ so that  the 
bubbling threshold is again respected. This implies that both $\phi$ and $\phi'$ are well-defined chain morphisms.

The next step is to show that the composition of the two chain morphisms $\phi'_{\eta}\circ \phi_{\eta}$ is chain homotopic to the identity as long as $\eta\in [\eta_{1},\eta_{2}]$.
The usual PSS technique applies to prove this statement. Again the only point worth making explicit concerns
the energy estimates. To discuss this, recall that the construction of the chain
homotopy between the identity and $\phi'_{\eta'}\circ \phi_{\eta}$ appeals to a new type of configuration
that we denote by $(r,\gamma, u, \gamma_1)$.
Here
$u:\R\to \widetilde{\mathcal{P}}_{0}(L_{\epsilon},L_{\epsilon})$,
$\gamma:(-\infty,0]\to \widetilde{L}_{\epsilon}$, $\gamma_{1}:[0,\infty)\to \widetilde{L}_{\epsilon}$ 
and with the notation $u'=p(u)$, $\gamma'=p(\gamma)$, $\gamma'_{1}=p(\gamma_1)$
 we have:
$$u'(\R\times\{0,1\})\subset L\ , \ \partial_{s}(u')+ \bar{J}(u')\partial_{t}(u')+\theta_{r}(s)\nabla
F(u')=0$$
$$
\frac{d\gamma'}{dt}=-\nabla \sigma(\gamma')\ ,\  \gamma(-\infty)=x \ , \
\gamma(0)=u(-\infty)~.~$$
$$
\frac{d\gamma_{1}'}{dt}=-\nabla \sigma(\gamma_{1}')\ ,\  \gamma_{1}(+\infty)=x' \ , \
\gamma_{1}(0)=u(+\infty)~.~$$
where $x$ and $x'$ are generators of $CW(\sigma,\widetilde{L}_{\epsilon};0)$.
The family of functions $\theta_{r}:\R\to [0,1]$ is  chosen so that when $r\to 0$ the family goes uniformly to $0$
and for sufficiently large $r$ it has support inside $[-r-1,r+1]$ and it is constant equal to $1$ in $[-r,r]$ and
is increasing in the interval $[-r-1,-r]$ and decreasing in the interval $[r,r+1]$.
It is again easy to estimate the energy of such configurations using the same formula as before. The conclusion in this case
is that because $x,x'$ are of weight $0$, we obtain $E(r,\gamma, u,\gamma_{1})\leq \osc(F) \leq \osc(H)+2\zeta < \delta'(V,J)$
so that we can deduce that $\phi'\circ\phi$ are chain homotopic to the identity by the usual PSS reasoning. 

Finally, we need to notice that $\bar{\psi}_{2,1}\circ \phi_{\eta_{2}}$ is chain homotopic to $\phi_{\eta_{1}}$. For this notice that the two Hamiltonians involved here are $F_{1}=\eta_{1}+\bar{\kappa}= -\osc(H)-\delta''+\bar{\kappa}$ and $F_{2}=\eta_{2}+\bar{\kappa}=\bar{\kappa}$. 
Thus $F_{1}$ and $F_{2}$ only differ by a constant. In particular, $F_{1}$ and $F_{2}$ have the same Hamiltonian flows. It follows that the difference between $\phi_{\eta_{1}}$
and $\phi_{\eta_{2}}$ only consists in the way the truncation is applied. In other words the actual
underlying moduli spaces are the same but when the respective chain morphisms are defined the truncations
take into account the difference between $F_{1}$ and $F_{2}$. But this is precisely the effect of 
$\bar{\psi}_{2,1}$. Indeed, up to chain homotopy, $\bar{\psi}_{2,1}$ is 
induced by a montone homotopy $G'$ relating $F_{2}$ to $F_{1}$ that has the form
 $G':[0,1]\times M\to \R$, $G'(t,x)= t\eta_{2}+(1-t)\eta_{1}+\bar{\kappa}(x)$. This implies that
 the comparison map is the identity at the level of the generators of the complexes $CF(L_{\epsilon},L_{\epsilon}; F_{i}, \bar{J})$ and its effect on the truncated complexes $CF^{a}_{b}$ just reflects the different truncation due to the difference between the action functionals $\mathcal{A}_{F_{2}}$ and $\mathcal{A}_{F_{1}}$.

\

This concludes
the proof of Proposition \ref{prop:shadow-est} in the case of simple cobordisms and under the simplifying 
geometric assumption in \S\ref{subsubsec:reposition}.

\

\

\

\subsection{Dropping the special assumptions.}\label{subsec:gen-case}
The statement of the proposition has been proved in the preceding sections under two 
assumptions:  that the cobordism $V$ is simple and that the projection of $V$ in the plane
is as described in  \S\ref{subsubsec:reposition} - Figure \ref{fig:elephant-in-boa}.

We will first see that the same method of proof applies in the case that $V$ is still simple but 
is not positioned as in  Figure \ref{fig:elephant-in-boa}.   To start this argument, first notice that in the 
proof discussed in the previous sections we can take without any difficulty
the points $Q, P\in \R+i$ so that $Re(Q)<-a$, $Re(P)>a$ for a constant
 $a>0$ as large as desired instead of assuming $Re(Q)\in (-2,-1)$, $Re(P)=\frac{3}{2}$. 
 
Assume now that $V$ is a simple cobordism. It is a simple exercise to show that there exists some constant 
$a>0$ and a symplectic diffeomorphism $\psi:\C\to \C$ with support in $[-a,a]\times \R$ 
so that if we put $\bar{\psi}=\psi\times id : \C\times M\to \C\times M$, then
$\tilde{V}=\bar{\psi}(V)$ has the projection $\pi(\tilde{V})$ as in Figure \ref{fig:elephant-in-boa}. 
We can now construct a Hamiltonian $\tilde{H}$, and $\tilde{V}'$, $\tilde{V}'',$ $\tilde{W}$  as described
in \S\ref{subsubsec:reposition} and  in \S\ref{subsec:comp-Lagr} (where the notation skips the $\tilde{\ }$)
but now starting from the cobordism $\tilde{V}$ and relative to two points $Q$ and $P$ so that
$Re(Q)<-a$, $Re(P)>a$. Next define $H=\tilde{H}\circ \bar{\psi}$,
$V'=\bar{\psi}^{-1}(\tilde{V}')$, $V''=\bar{\psi}^{-1}(\tilde{V}'')$, $W=\bar{\psi}^{-1}(\tilde{W})$.
Notice that $\mathcal{S}(\tilde{V})=\mathcal{S}(V)$ so that the oscillation of $H$ continues to be 
controlled by the shadow of $V$. From this point on we continue the proof exactly as in the sections
\S\ref{subsec:comp-Lagr} - \S\ref{subsubsec:PSS} but for  $H$, $V'$, $V''$, $W$ as just defined.  
Notice, in particular, that because the points $Q$ and $P$ are away from the support of $\psi$ the arguments involving the behaviour of holomorphic curves near $Q$ and $P$ do not need any adjustment. 
This ends the proof of the proposition in the case of a simple cobordism.

The case of a non-simple cobordism is basically
 a trivial consequence of the proof in the simple case because the condition $L\cap L_{i}=\emptyset$
 $\forall \ 1\leq i\leq k$ and $L\cap L'_{j}=\emptyset$ $\forall \ 1\leq j\leq s$ insures that the ends
 of the cobordism $V$ that are different from both $L$ and $L'$ do not interfere in any way in the proof.
 This concludes the proof of Proposition \ref{prop:shadow-est}.
 
\section{Proof of Theorem \ref{thm:metrics}}\label{sec:proofs}
 Each point of the Theorem is the subject of one of the  subsections below.

\subsection{The metrics $d^{\ast}$ for $\ast\geq w-m$} 
We first remark that formula (\ref{eq:distances}) defines a pseudo-metric. To start, notice that cobordisms can 
be glued.  More precisely, if $V: L\cobto L'$, $V':L'\cobto L''$ are two cobordisms, assumed both to be cylindrical outside 
$[0,1]\times \R$, then $V'':L\cobto L''$ is 
obtained as the union $\bar{V}\cup \tilde{V'}$ where $\bar{V}= V\setminus (-\infty, -2]\times \R\times M$
and $\tilde{V'}=\{(x-3, y, m)  \ | \ (x,y, m)\in  V'\setminus [2,+\infty)\times \R\times M\}$. 
It is easy to see that $\mathcal{S}(V'')=\mathcal{S}(V)+\mathcal{S}(V')$. As a consequence $d^{\ast}$ verifies the triangle inequality. Further, for each cobordism $V: L \cobto L'$ we can use the planar transformation
$z \to - z$ followed by an appropriate translation in $\C$ to construct a cobordism $\hat{V}:L'\cobto L$ so that
$\mathcal{S}(V')=\mathcal{S}(V)$. As a consequence $d^{\ast}$ is symmetric and thus a pseudo-metric.

We now need to see that when $\ast\geq w-m$ this pseudo-metric is non-degenerate.  Consider two distinct Lagrangians
$L,L'\in \mathcal{L}^{\ast}(M)$ and a cobordism $V:L\cobto L'$, $V\in \mathcal{L}^{\ast}_{cob}(\C\times M).$ By Proposition \ref{prop:shadow-est}, for arbitrarily small $\epsilon_{0}$, there exists an almost complex structure $J$ such that $d^{\ast}(L,L')\geq \min\{ w(L,L') - \epsilon_{0}, \delta(V; J)\}$. Given that $L\not=L'$, $w(L,L')>0$.
For  $\ast=w-m$, if $u$ is a $J$-holomorphic disk or sphere (as in (\ref{eq:threshold})), then $\omega(u) = \rho \mu(u)$ and as $\mu(u)\in \Z$ it follows that $\omega(u)\geq |\rho |$.  Thus $\delta(V;J)$ is bounded 
from below by $|\rho|$ (and is equal to $+\infty$ when $\rho = 0$). This argument also applies for 
$\ast=m$ as well as for $\ast=e$. Thus if $\ast=w-m,m,e$ we obtain that the pseudo-metric $d^{\ast}$
is non-degenerate. The Hamiltonian orbit case  reduces to a statement 
that is well-known but we will sketch a direct proof for completeness in the next section.

\subsection{Relation to the Hofer norm.}\label{subsubsec:rel-H-norm}

Let $\phi\in \Ham(M,\omega)$.  The Hofer norm of $\phi$ is given by:
$$ ||	\phi||_{H}=\ \ \inf_{G,\ \phi^{G}_{1}=\phi}\ \  \int_{0}^{1}(\max_{x\in M} G(t,x)-\min_{x\in M} G(t,x)\ )dt$$
where $G:[0,1]\times M\to \R$ is a Hamiltonian function. Here $M$ is assumed to be compact or $\phi$ and 
$G$ are assumed to have compact support.  There is also a variant of the associated distance that is relative to 
a closed Lagrangian $L\subset M$. Assume that $L'$ is a Lagrangian that is Hamiltonian isotopic to $L$.
Then we can define a pseudo-distance:
$$d_{H}(L,L')=\ \ \inf_{G,\ \phi^{G}_{1}(L)=L'}[\int_{0}^{1}(\ \max_{z\in \phi^{G}_{t}(L)} G(t,z) -\min_{_{z\in \phi^{G}_{t}(L)}}G(t,z)\ )dt]~.~$$
This pseudo-distance is actually a distance on the Hamiltonian orbit of $L$ - the Hofer distance between Hamiltonian-isotopic Lagrangians - that has been considered by many authors starting from \cite{ChekanovFinsler} in slightly different forms (herein we follow \cite{Bar-Cor:Serre,Bar-Cor:NATO,Cha:thesis}; see \cite{UsherObservations} for additional references). 
To show that $d_{H}$ is non-degenerate one can proceed just as in our  argument for the first point of Proposition \ref{prop:shadow-est}. For completeness, we sketch the proof below.

We want first to show that if $L=\phi^{G}_{1}(L')$, then for $\epsilon_0 >0$ there exists an almost complex structure $J$ on $M$ such that
 $$\int_{0}^{1}(\max_{x\in M} G(t,x)-\min_{x\in M} G(t,x)\ )dt \geq \min \{\delta(L, J),  w(L,L')-\epsilon_0\}~.~$$
We first  pick $J$ as $J_{-}$ was chosen in \S\ref{subsec:comp-Lagr}. We assume that the quantity $$\osc(G)=
 \int_{0}^{1}(\max_{x\in M} G(t,x)-\min_{x\in M} G(t,x))dt$$  is smaller than the bubbling threshold $\delta(L,J)$.
 By a naturality transformation such as $\Psi'$ in (\ref{eq:nat2}) together with the module action from \S\ref{subsec:Fl-hlgy-Morse}, the proof reduces to show that 
 \begin{equation}\label{eq:est-mod11}
 R\ast HF^{a}_{b}(L,L; G, \bar{J})\not=0
 \end{equation} where $(\phi^{G}_{t})_{\ast}\bar{J}=J$.
 This non-vanishing is in perfect analogy to formula (\ref{eq:mod-est10}). Here $a= \int_{0}^{1}(\max_{x\in M} G(t,x))dt +\zeta$ and $b=\int_{0}^{1}(\min_{x\in M} G(t,x)\ )dt-\zeta$ where $\zeta$ is an arbitrarily small constant.  To show (\ref{eq:est-mod11})  we apply the construction of the PSS maps as described in \S\ref{subsubsec:PSS}. 
 The energy estimates in this case  are exactly what is required for the argument to work. Indeed, in the argument in 
 \S\ref{subsubsec:PSS} we dealt with a time independent Hamiltonian $F$ but if we apply the same energy calculations to a time dependent Hamiltonian $G$, 
 then the oscillation of $F$ is replaced by $\osc(G)$, $\sup(F)$ by $\int_{0}^{1}(\max_{x\in M} G(t,x))dt$ and $\inf(F)$ by $\int_{0}^{1}(\min_{x\in M} G(t,x)\ )dt$. The final step is to see that  $\max$ and $\min$ in these formulas
 can be taken as in the definition of $d_{H}$, that is over $\phi^{G}_{t}(L)$. But this is easy to do by 
 truncating $G(t,z)$ for each $t$ outside a neighbourhood of $\phi^{G}_{t}(L)$.

\

The shadow of  cobordisms and the Hofer norm are naturally related through the  Lagrangian suspension 
construction.  

\

Fix a connected, closed Lagrangian $L\subset M$ and let its Lagrangian suspension along $G$ be $L^{G}$ as in 
\S\ref{subsubsec:orb}. It is immediate to see from
Definition \ref{def:shadow} that:
\begin{equation}\label{eq:Lag-Susp}
\mathcal{S}(L^{G})=\int_{0}^{1}[\ \max_{z\in \phi^{G}_{t}(L)} G(t,z) -\min_{_{z\in \phi^{G}_{t}(L)}}G(t,z)\ ]dt~.~
\end{equation}

Thus, for a Lagrangian $L_{0}$ in $M$, the metric $d^{L_{0}}$ defined on the orbit
$\mathcal{L}^{L_{0}}(M)$ of $L_{0}$ under the action of the Hamiltonian group, as provided by Theorem \ref{thm:metrics}, satisfies:
 $$d^{L_{0}}(L,L')= d_{H}(L,L')$$
 for any $L,L'\in \mathcal{L}^{L_{0}}(M)$. 

\begin{rem}
 For any non-constant $\phi\in \Ham(M,\omega)$ there is some Lagrangian $L$ in $M$ - possibly taken in a sufficiently small Darboux chart -  so that $\phi(L)\not= L'$. Thus, the non-degeneracy of the Lagrangian metric $d_{H}(-, -)$
implies that  the Hofer norm on $\Ham(M,\omega)$ itself is non-degenerate. This method to show the non-degeneracy of the Hofer norm
 is due to Polterovich \cite{Po:LagrDispl} and Chekanov \cite{Chekanov:Lag-energy-1}.
\end{rem}

\subsection{Surgery and non-isotopic Lagrangians at finite $d^{w-m}$ distance.} 
An example of two non-isotopic Lagrangians and a weakly-monotone cobordism relating them was
constructed in \cite{Bi-Co:cob1}. We give here an outline of the construction because we want to also discuss
the shadow of this cobordism.

We start by recalling the Lagrangian surgery construction, \cite{La-Si:Lag}, \cite{Po:surgery}. This is based on a simple local construction.  Fix the following two Lagrangians:
$L_{1}=\R^{n}\subset \C^{n}$ and $L_{2}= i\R^{n}\subset \C^{n}$.

\begin{figure}[htbp]
   \begin{center}
      \epsfig{file=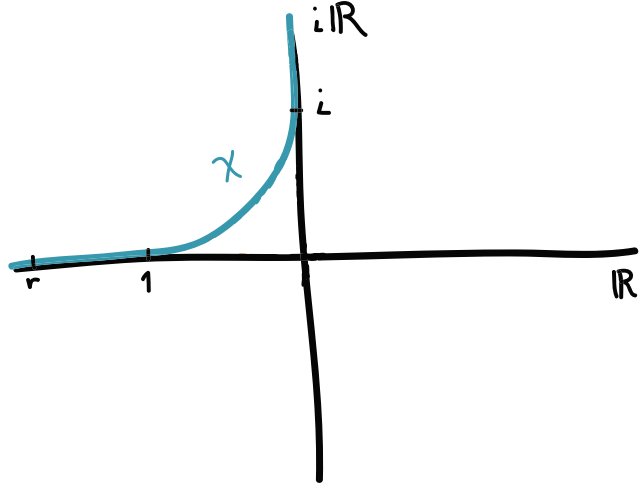, width=0.50\linewidth}
   \end{center}
   \caption{\label{fig:handle} The curve $\chi\subset \C$.}
\end{figure}

Consider a curve  $\chi\subset \C$, $\chi(t)=a(t)+ ib(t)$,
$t\in\R$, so that (see also Figure
\ref{fig:handle}):
$\chi$ is smooth; $(a(t), b(t))=(t,0)$ for $t\in(-\infty,-1]$; $(a(t), b(t))=(0, t)$ for $t\in [1,+\infty)$; 
$a'(t), b'(t)>0$ for $t\in (-1,1)$.
 Let $$L =
\Bigl\{\bigl((a(t)+ib(t))x_{1}, \ldots, (a(t)+ib(t))x_{n}\bigr) \mid t
\in \mathbb{R}, \, \sum x_{i}^{2}=1\Bigr\}\subset \C^{n}~.~$$
 It is easy to see that $L$ is Lagrangian.
 By an abuse of notation because we omit the handle $\chi$ from the
notation and  we will denote $L=L_{1}\# L_{2}$. Notice that different choices of 
handles $\chi$ produce Hamiltonian isotopic Lagrangians $L$ (for $n>1$). By choosing the handle small enough, we can have the result of the surgery be contained in an arbitrarily small neighbourhood of $L_{1}\cup L_{2}$.

\begin{figure}[htbp]
     \begin{center}
        \epsfig{file=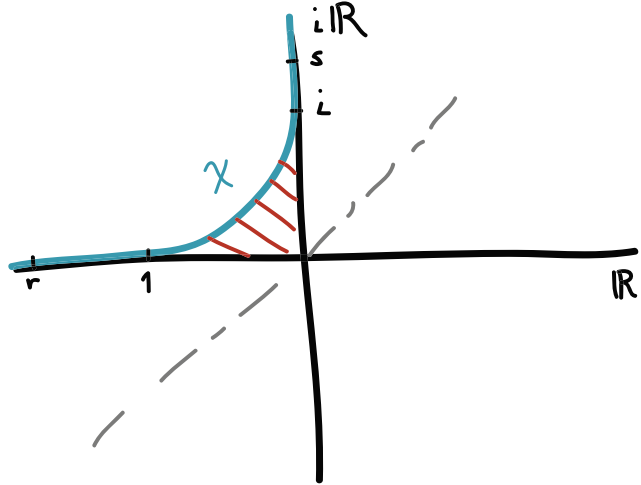, width=0.50\linewidth}
     \end{center}
     \caption{\label{fig:projection-image} The projection of $W_{0}$ is
       the red region together with the two semi-axes
       $(-\infty,0]\subset \R$ and $i[0,+\infty)\subset i\R$ and the
       curve $\chi$.}  
  \end{figure}

\begin{figure}[htbp]
     \begin{center}
        \epsfig{file=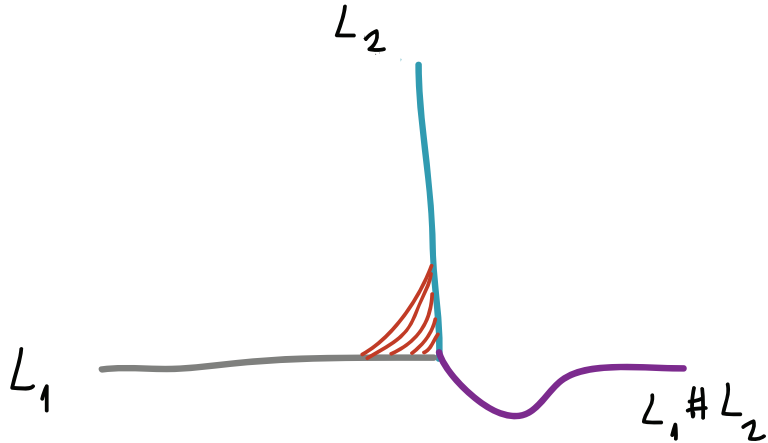, width=0.60\linewidth}
     \end{center}
     \caption{\label{fig:surgery-trace} The trace of the surgery after
       projection on the plane.}
  \end{figure}

 There is a Lagrangian cobordism $L \cobto
   (L_{1},L_{2})$ constructed as follows.
   Define 
   $$\widehat{\chi}=\Bigl\{ \bigl( (a(t)+ib(t))x_{1}, \ldots,
   (a(t)+ib(t))x_{n+1} \bigr) \mid t\in \R, \, \sum x_{i}^{2}=1
   \Bigr\} \subset \C^{n+1}$$ and notice that 
   $\widehat{\chi}$ is also Lagrangian.
 Consider the projection {$\pi: \C^{n+1}\to \C$, $\pi(z_{1},\ldots
   z_{n+1})=z_{1}$} and denote by $\widehat{\pi}$ its restriction
   to $\widehat{\chi}$.
    Define $W =
   \widehat{\pi}^{-1}(S_{+})$ where $S_{+}=\{(x,y)\in \R^{2} \ | \
   y\geq x\}$, see Figure \ref{fig:projection-image}. (As usual, we
   identify $\R^{2}$ with $\C$ under $(x,y)\to x+iy$.)  
  Consider $W_{0}=W\cap \pi^{-1}([-2,0]\times[0,2])$. It is
  not difficult to see that $W_{0}$ is a manifold with boundary and that 
  $\partial W_{0}=\{(-2,0)\}\times L_{1}\cup \{ (0,2)\}\times L_{2} \cup 
  \{0,0\}\times L$.  To finish the construction of the cobordism we adjust $W_{0}$ (as described explicitly in \cite{Bi-Co:cob1}) so as to continue the $L$-boundary component to be cylindrical.
  The resulting Lagrangian $W'$ provides the cobordism desired between $L$ and
  $(L_{1},L_{2})$ - see also Figure \ref{fig:surgery-trace}.

Going from the local argument above to a global one is easy. Suppose
that we have two Lagrangians $L'$ and $L''$ that intersect
transversely, possibly in more than a single point. At  each
intersection point we fix symplectic coordinates mapping (locally)
$L'$ to $\R^{n}\subset \C^{n}$ and mapping (again locally) $L''$ to
$i\R^{n}\subset \C^{n}$.  We then apply the construction above at each
of these intersection points.  This produces a new Lagrangian
submanifold $L'\#L''$ as well as a cobordism
$W'':L'\#L''\cobto (L',L'')$ (caveat: $L'\#L''$ is topologically not a connected sum if there are
several intersection points).  

\

We are interested in the shadows of the Lagrangians
resulting from this construction. There are two useful remarks in this direction. 
First, for the cobordism $W''$ above, we see easily that for any $\epsilon>0$
we can find a sufficiently small handle $\chi$ so that $\mathcal{S}(W'')\leq \epsilon$. However,
different handles lead to different outputs of the surgery as the resulting $L'\# L''$ are Hamiltonian
isotopic (for $n>1$) but not identical. Secondly, if in the place of  $L'$ and $L''$ we take cobordisms $V'=\gamma'\times N'$,
$V''=\gamma''\times N''$ where $\gamma'$ and $\gamma''$ are appropriate curves in $\C$ that intersect
transversely at a single point, then for any $\epsilon$ we may construct the surgered Lagrangian $V'\# V''$ 
so that its shadow is smaller than $\epsilon$ by again using in the construction a sufficiently small handle $H$. 

\

Our example is based on the construction decribed above. We will start our construction in the ambient manifold $M'=\C\backslash\{P_{1},P_{2},P_{3}\}$ where $P_{1}$, $P_{2}$, $P_{3}$ are three points $P_{i}\in \C$. One may worry that excising points creates concave ends in our symplectic manifold, but this is easily overcome by a positivity of intersections argument. Alternatively one may glue in small handles at the three points $P_{1}$, $P_{2}$, $P_{3}$. Finally, one can work inside a large ball, and compactify it to a sphere, adding a handle at infinity. This gives us an example in the closed surface of genus $4.$ The following arguments apply uniformly in all cases.

\begin{figure}[htbp]
   \begin{center}
      \epsfig{file=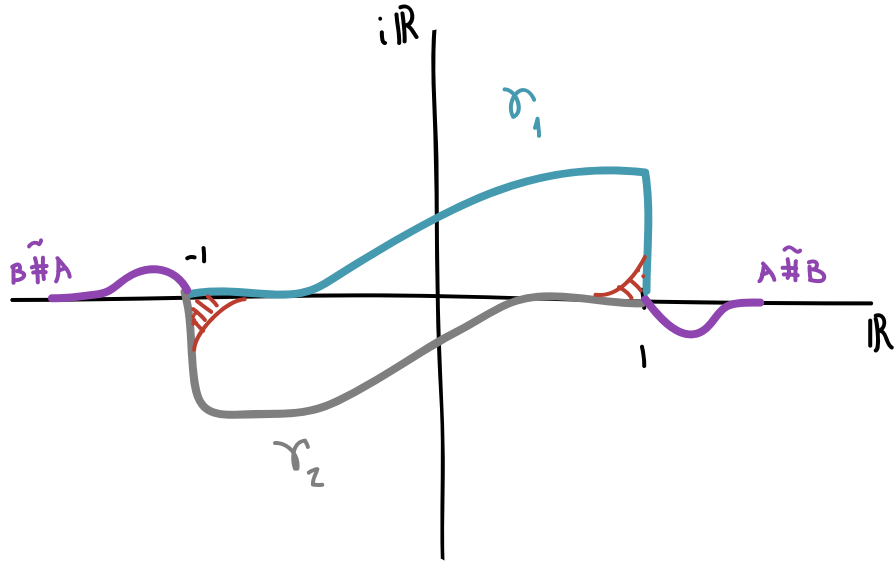, width=0.70\linewidth}
   \end{center}
   \caption{\label{fig:plane-surgery1} The projection of $V$ on $\C$;
     in {\color{red} red} the surgery regions; the curves $\gamma_{1}$
     (in blue) and $\gamma_{2}$ (in gray).}
\end{figure}

\begin{figure}[htbp]
   \begin{center}
      \epsfig{file=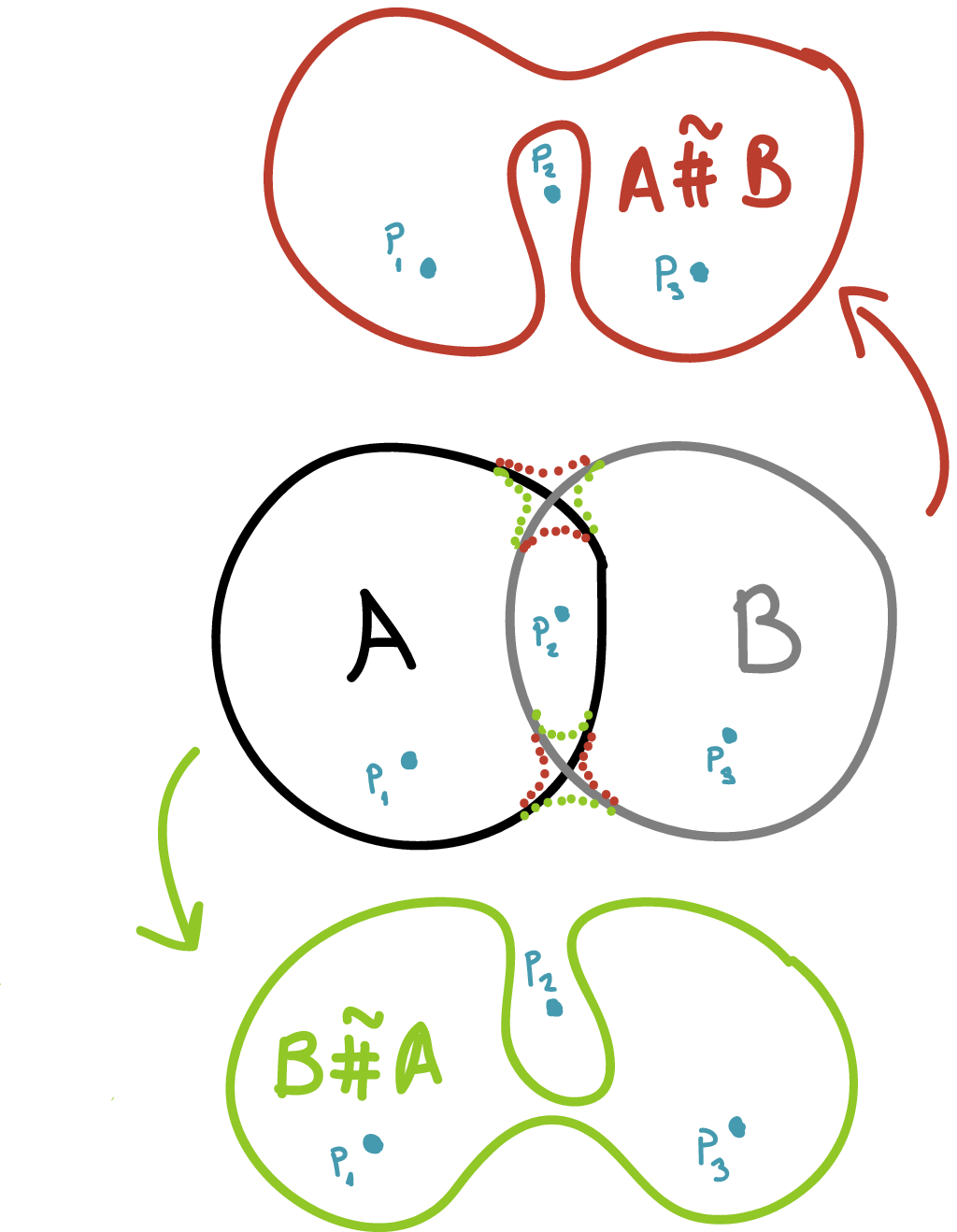, width=0.40 \linewidth}
   \end{center}
   \caption{\label{fig:surgered-circles} The two circles $A$ and $B$
     as well as {\color{red} $A\# B$} and {\color{green}
       $B\# A$}. The three puncture points are indicated in
     blue.}
\end{figure}

We consider two circles $A=\{z\in \C \ :\ |z+1/2|=1\}$ and $B=\{x\in \C\
: \ |z-1/2|=1\}$ and denote by $D(A)$ and $D(B)$  the two
disks bounded by $A$ and $B$ respectively. We assume that the positions of the points $P_{i}$ relative to the circles $A,B$ are such that $P_{1}\in D(A)\backslash D(B)$, $P_{2}\in
D(A)\cap D(B)$ and $P_{3}\in D(B)\backslash D(A)$ as at the middle of 
Figure~\ref{fig:surgered-circles}.    We consider two smooth curves in
the plane $\C$, $\gamma_{1}:[-1,1]\to \C$ and $\gamma_{2}: [-1,1]\to
\C$ so that - see Figure \ref{fig:plane-surgery1}:
\begin{itemize}
  \item[i.] $\gamma_{1}(t)=t$ for $t\in [-1,-1/2]$
  \item[ii.] $\gamma_{1}(t)=1+(1-t)i$ for $t\in [1/2,1]$
  \item[iii.] $Re(\gamma_{1}(t))$ is strictly increasing for $t\in
   (-1/2,1/2-\epsilon)$. $Im(\gamma_{1}(t))$  is strictly increasing for $t\in
   (-1/2,1/2-\epsilon)$ and strictly decreasing for $t\in
   (1/2-\epsilon, 1/2)$.
  \item[iv.] $\gamma_{2}(t)=-\gamma_{1}(t)$ for all $t\in [-1,1]$.
\end{itemize}

We now consider the Lagrangians $A'=\gamma_{2}\times A \subset
\mathbb{C} \times M$ and $B'=\gamma_{1}\times B \subset \C\times M$.
By performing surgery
at both intersection points $A\cap B$ we can extend the union of the
two Lagrangians $A'\cup B'$ towards the positive end as well as
towards the negative end as in the Figure \ref{fig:plane-surgery1}
thus obtaining a cobordism $V:A\# B\cobto B\# A$.

Put $L=A\# B$ and $L'=B\# A$.  It is easy to see that $L$ and $L'$ look as in
Figure~\ref{fig:surgered-circles} and are exact.
We notice that $L$ and $L'$ are not smoothly isotopic in $M'$. 
However, it is shown in \cite{Bi-Co:cob1} that, by choosing the handles associated with the surgeries in the two intersection points of $A$ and $B$ appropriately, $V$ can be made monotone with minimal Maslov class $1$. 
It is also clear that the shadow of $V$ can be made as close as desired to the sum of two areas,
one corresponding to the shadow of one handle used at the $A\# B$ end of the surgery and the other
corresponding to the shadow of the handle used at the $B\# A$ end  -  as suggested by Figure \ref{fig:plane-surgery1}. Notice however that diminishing the size of these handles also modifies  the 
ends $L=A\# B$ and $L'= B\# A$. Indeed, Proposition \ref{prop:shadow-est} 
shows that if $L$ and $L'$ are fixed the shadow of the cobordism $V$ can not be arbitrarily small.

It is also noticed in \cite{Bi-Co:cob1} that  $L$ and $L'$ are not cobordant via a monotone cobordism $V$ with $N_{V}\geq 2$. In short, this is shown by observing that  the Floer homologies $HF(S ,L)$ and $HF(S,L')$
 where $S$ is the vertical semiaxis in $\C$, pointing up and starting at $P_{2}$ are non-isomorphic. But by  the
 results in \cite{Bi-Co:cob1,Bi-Co:cob2} if a monotone simple cobordism would relate $L$ and $L'$ then 
 for any other Lagrnagian $N$ we would have  $HF(N, L)\cong HF(N,L')$ (indeed, $L$ and $L'$ would even be isomorphic in the appropriate Fukaya category).
  
  \
  
 In summary, the two exact Lagrangians $L,L'\in \mathcal{L}^{e}(M')$ constructed before
  are not smoothly isotopic and 
  $$0<d^{w-m}(L,L')<\infty\ , \ d^{m}(L,L')=\infty~.~$$
  We will see in \S\ref{subsubec:deg} that we also have $d^{g}(L,L')=0$.
  
\subsection{The pseudo-metric $d^{g}$ is degenerate.}\label{subsubec:deg}

We consider here the pesudo-metric $d^{g}$ given as in the formula in Theorem \ref{thm:metrics} but for $\ast=g$. Recall that $g$ stands for general. In other words, there is no constraint imposed on either the Lagrangians or the cobordisms involved.

\ 

The following construction was suggested to us by Emmy Murphy.

\

Consider two Lagrangians $L$ and $L'$ and a cobordism $V:L\cobto L'$.

For $\epsilon >0$ define a map $\tilde{\epsilon}:\C \times M \to \C\times M$ given by rescaling in the imaginary direction in the plane $\tilde{\epsilon}(x+iy,z) = (x + i \epsilon y, z).$ Put $V^{\epsilon} = \tilde{\epsilon}(V).$ Denote by $e: V \to \C \times M$ the natural embedding, and define a smooth embedding $e_{\epsilon}:V \to \C \times M$ by $e_{\epsilon} = \tilde{\epsilon} \circ e.$ 

Let $R$ be a compact rectangle in $\C$ outside whose interior $V$ is cylindrical: for example one can take $R = (\beta_{-},\beta_{+}) \times (Y_{-},Y_{+}),$ where $Y_{-} < \inf \mathrm{Im}(\pi(V) \cap [\beta_{-},\beta_{+}]),$ and $Y_{+} > \sup \mathrm{Im}(\pi(V) \cap [\beta_{-},\beta_{+}])$. Clearly, $V^{1}=V$ and $$\mathrm{Area} (ou_{V^{\epsilon}})=\epsilon\ \mathrm{Area}(ou_{V}) < 
\epsilon \mathrm{Area}(R).$$  

For a point $v \in \C \times M,$ consider the symplectomorphism $l_{\epsilon,v}: \C \times M \to \C \times M$ given by $$l_{\epsilon, v}(x+iy,z)=(x+iy, z) - (i(1-\epsilon)\mathrm{Im}(v), 0)~.~$$ Viewed as a family of maps, $l_{\epsilon, v}$ depends smoothly on both parameters.
Note that $l_{\epsilon,v}(v) = \tilde{\epsilon}(v).$  Hence, $dl_{\epsilon,v}: T_{v}(\C\times M) \to T_{\tilde{\epsilon}(v)}(\C\times M)$ is an isomorphism of symplectic vector spaces.

Define a smooth bundle map
$\phi^{\epsilon}: TV \to e_{\epsilon}^{\ast} T(\C\times M)$ by
$$\phi^{\epsilon}(\xi)=dl_{\epsilon,e(a)}\circ de (\xi)$$ for each
$\xi\in T_{a}V$. The image of this bundle map is a Lagrangian sub-bundle
and, because this bundle map covers the smooth embedding $e_{\epsilon}$,
 we may apply the Gromov-Lees h-principle (\cite{Gr-hprinciple,Lee-h}) to $e_{\epsilon}$. As a result we obtain that  
 $e_{\epsilon}$ can be approximated arbitrarily well in $C^{0}$ norm by Lagrangian immersions.
 In particular, for any $\delta>0$ we may find a Lagrangian immersion
 $e_{\epsilon,\delta}: [0,1]\times L\to \C\times M$ with image $V^{\epsilon, \delta}$ 
 so that $\mathrm{Area}(ou_{V^{\epsilon,\delta}})\leq \epsilon \mathrm{Area}(R) +\delta$.

We now modify $V^{\epsilon,\delta}$ twice: first we perturb the immersion to a new immersion with only 
transverse double points and, secondly, surger all the self-intersection points by using very small 
Lagrangian handles so as to get an embedded Lagrangian $V^{\epsilon,\delta,\delta'}$
so that $\mathrm{Area}(ou_{V^{\epsilon,\delta,\delta'}})\leq \epsilon \mathrm{Area}(R)+\delta +\delta'$. 

By taking the (generic) perturbation of the immersion $e_{\epsilon,\delta}$ small enough and
by taking the surgery handles to be also small enough, we may assume that $\delta'\leq \epsilon$ and $\delta\leq \epsilon.$ Hence we get a cobordism $V^{\epsilon,\delta,\delta'}:L\cobto L'$ with
$$\mathrm{Area}(ou_{V^{\epsilon,\delta,\delta'}})\leq  \epsilon\ (\mathrm{Area}(R) + 2)~.~$$
Therefore $d^{g}(L,L')=0$ and thus the pseudo-metric $d^{g}$ is degenerate.

\

This concludes the proof of Theorem \ref{thm:metrics}. 
\section{Additional comments}\label{sec:comm}

\subsection{Relation to spectral distance.} 
The argument for the proof of Proposition \ref{prop:shadow-est} suggests that in the setting where $L$ and $L'$ are Hamiltonian isotopic and exact (the same would hold in the weakly exact case: $\omega, \mu|_{\pi_{2}(M,L)}=0$), assuming that the cobordism is monotone, one can replace $w(L,L')$ in the statement of the Proposition by $d_{S}(L,L')$, the spectral distance between $L$ and $L'$ (introduced in \cite{Vi:generating-functions-1}, see also \cite{HLS-CoisotrRigidity} for additional references). For a fixed Lagrangian $L$ recall from the work of Milinkovic \cite{Mil} that,  if $L'$ is sufficiently $C^1$-close to $L$, then  $d_{S}(L,L')=d_{H}(L,L')$. Therefore, we expect that, at least under this additional proximity assumption, $d^{\ast}(L,L')=d_{H}(L,L')$  for all $\ast\geq m$. 

\subsection{Lower bound for the shadow in the monotone case.} We believe that an adaptation of the proof of
Proposition \ref{prop:shadow-est} shows that, under the assumptions of the Proposition, and if, additionally, $ L, L'\in \mathcal{L}^{m}(M)$, $V\in \mathcal{L}^{m}_{cob}(\C\times M),$ then we have:
$$\mathcal{S}(V)\geq w(L,L')~.~$$
This inequality fits with the leitmotiv of the paper:  more rigid topological constraints lead to
sharper inequalities. Indeed, in the setting of the proposition, this expected inequality shows
that monotonicity is sufficient to eliminate $\delta(V,J)$ from the general inequality given by Proposition \ref{prop:shadow-est}. Further, as seen in \S\ref{subsubsec:rel-H-norm}, if we assume that $V:L\cobto L'$ is a Lagrangian suspension, then the inequality becomes even stronger as we can replace $w(L,L')$ by the Hofer distance $d_{H}(L,L')$ between $L$ and $L'$.

\subsection{Categorical view-point}  A somewhat more conceptual perspective on the construction of the metrics $d^{\ast}$ from Theorem
\ref{thm:metrics} is as follows.
Using the notion of cobordism, one can define - as in \cite{Bi-Co:cob1} - various categories 
that have as objects Lagrangians in $M$ and have morphisms given by Lagrangian cobordisms.
As before, the specific Lagrangians and cobordisms involved are subject to the constraints encoded in the superscript $-^{\ast}$, where $\ast$ can be any of the 
conditions listed in \S\ref{subsec:class-Lagr}.
The simplest such category, $\mathcal{L}ag^{\ast}_{s}(M)$, has
as objects the Lagrangians in
$\mathcal{L}^{\ast}(M)$ and as morphisms the horizontal isotopy classes of simple cobordisms 
$V:L\cobto L'$ so that  $V\in \mathcal{L}^{\ast}_{cob}(C\to M)$ .

Given a small category $\mathcal{C}$ assume that the morphisms of 
$\mathcal{C}$ are endowed with a valuation $\nu: \mor_{\mathcal{C}}\to [0,\infty)$
in the sense that $\nu (f\circ g)\leq \nu (f)+\nu (g)$ for all composable morphisms $f$ and $g$,
 and, for each morphism $f\in \mor (A,B)$, there
 is a morphism $\bar{f}\in \mor (B,A)$ with $\nu(f)=\nu(\bar{f})$. Such a valuation induces a 
 pseudo-metric metric $d_{\nu}$ on $\mathcal{O}b(\mathcal{C})$ that is given by:
 $$d_{\nu}(A,B)=\inf_{\varphi \in \mor (A,b)}\nu (\varphi)~.~$$
 This number is taken to be infinite in case there are no morphisms from $A$ to $B$.
 In case the valuation is non-degenerate in the sense that $\nu(f)=0$ iff $f=id_{X}$ for some object $X$,
 then the pseudo-metric is a true metric (with this definition the metric is finite only for objects that are related by some morphism).

The shadow of cobordisms, as given in Definition \ref{def:shadow},
provides a valuation on the category $\mathcal{L}ag^{\ast}_{s}(M)$ by putting for each
morphism $[V]$ represented by a cobordism $V$ :
\begin{equation} \label{eq:size-mor}
\nu([V])=\inf_{V'}\ \{\mathcal{S}(V') : V'\ \mathrm{horizontally\ isotopic\ to\ } V\}~.~
\end{equation}
Obviously, Theorem \ref{thm:metrics} shows that the resulting pseudo-metrics $d^{\ast}=d_{\nu}^{\ast}$ are non-degenerate for $\ast\geq w-m$ and degenerate for $\ast=g$. 

\subsection{Immersed Lagrangian cobordism} Following the work of Akaho \cite{Aka}
as well as Akaho-Joyce \cite{Ak-Jo} (see also \cite{Alston-Bao}), Floer theory is also defined for a class of immersed Lagrangians
so-called unobstructed. The cobordism machinery can also be adapted without any trouble to this 
setting and we expect that there are variants of both Proposition \ref{prop:shadow-est} and Theorem \ref{thm:metrics} in this context.

\bibliography{bibliography}

\end{document}